\documentclass[onecolumn]{elsart3p}

\usepackage{graphicx}

\usepackage{amssymb,amsmath,bm}   

\usepackage{lineno}

\newcommand{\bx}{\textbf{x}}
\newcommand{\by}{\textbf{y}}
\newcommand{\bs}{\textbf{s}}

\newcommand{\bi}{\begin{itemize}}
\newcommand{\ei}{\end{itemize}}
\newcommand{\ben}{\begin{enumerate}}
\newcommand{\een}{\end{enumerate}}
\newcommand{\be}{\begin{equation}}
\newcommand{\ee}{\end{equation}}
\newcommand{\bea}{\begin{eqnarray}}
\newcommand{\eea}{\end{eqnarray}}
\newcommand{\bc}{\begin{center}}
\newcommand{\ec}{\end{center}}
\newcommand{\spl}[2]{\left\{\begin{array}{ll}#1\\#2\end{array}\right.}
\newcommand{\ie}{{\it i.e.\ }}
\newcommand{\eg}{{\it e.g.\ }}

\newcommand{\tbox}[1]{{\mbox{\tiny #1}}}
\newcommand{\mbf}[1]{{\bm #1}}           
\newcommand{\pO}{{\partial\Omega}}
\newcommand{\uN}{u^{(N)}}
\newcommand{\emach}{\epsilon_\tbox{mach}}
\newtheorem{rmk}[thm]{Remark}

\newtheorem{cnj}[thm]{Conjecture}
\newtheorem{lemma}[thm]{Lemma}
\newcommand{\bp}{{\bf Proof:\ }}   
\newcommand{\ep}{\hfill $\square$ \vspace{1ex} \\} 
\newcommand{\bmal}{{\bm \alpha}}             

\newcommand{\cu}{\tilde{u}}            
\newcommand{\cOmega}{\tilde{\Omega}}   
\newcommand{\dmax}{D_\tbox{max}}       

\begin{document}

\begin{frontmatter}


\title{Stability and convergence of the Method of Fundamental Solutions for
Helmholtz problems on analytic
domains}
\author[AHB]{A.~H.~Barnett\corauthref{co:AHB}}
\ead{ahb@math.dartmouth.edu}
\ead[url]{www.math.dartmouth.edu/$\sim$ahb}
\author[TB]{T.~Betcke}
\ead{timo.betcke@manchester.ac.uk}
\ead[url]{www.maths.man.ac.uk/$\sim$tbetcke}
\corauth[co:AHB]{Corresponding author. tel:+1-603-646-3178.
fax:+1-603-646-1312}


\address[AHB]{Department of Mathematics, 6188 Kemeny Hall, Dartmouth
  College, Hanover, NH, 03755, USA}
\address[TB]{School of Mathematics, The University of Manchester,
  Manchester, M13 9PL, UK}



\begin{abstract}
The Method of Fundamental Solutions (MFS) is a popular tool to solve Laplace
and Helmholtz boundary value problems. Its main drawback is that it often leads
to ill-conditioned systems of equations. In this paper we investigate for the
interior Helmholtz problem
on analytic domains how the singularities (charge points) of the MFS basis
functions have to be chosen such that approximate solutions can be represented
by the MFS basis in a numerically stable way. For Helmholtz problems on the
unit disc we give a full analysis
which includes the high frequency (short wavelength) limit.
For more difficult and nonconvex
domains such as crescents we demonstrate how the
right choice of charge points is connected to how far
into the complex plane the solution of the boundary value problem can be
analytically continued, which in turn depends on both domain shape
and boundary data.
Using this we develop a recipe for locating charge points which
allows us to reach error norms of typically $10^{-11}$ on a wide variety
of analytic domains.
At high frequencies of order only 3 points per wavelength are
needed, which compares very favorably to boundary integral methods.
\end{abstract}

\end{frontmatter}









\section{Introduction}
\label{sec:introduction}

The Method of Fundamental Solutions (MFS), also known as the charge
simulation method or the method of auxiliary sources,
is a well known method for solving Laplace or Helmholtz
boundary value problems (BVPs).
The idea is to approximate the solution by fundamental
solutions of the Laplace or Helmholtz equation whose singularities lie outside
the domain. Consider the boundary value problem
\begin{subequations}
\label{eq:bvproblem}
\begin{align}
\Delta u+k^2u&=0\quad\text{in }\Omega \label{eq:bvproblema},\\
 u&=v\quad\text{on }\partial \Omega\label{eq:bvproblemb},
\end{align}
\end{subequations}
where $\Omega\subset\mathbb{R}^2 = \mathbb{C}$
is a simply connected planar domain with
analytic boundary $\partial\Omega$.
Recall that
the solution is unique if and only if $k^2$ is not a Dirichlet
eigenvalue (of the Laplace operator) for the domain;
physically this is a resonance effect.
%
The idea of the MFS is to approximate $u$ by a
linear combination of fundamental solutions of the form
\begin{equation}
\label{eq:mfs}
u(\bx)\approx
\uN(\bx)=
\frac{i}{4}\sum_{j=1}^{N}\alpha_j H_0^{(1)}(k|\bx-\by_j|),\quad   
\by_j\in\mathbb{R}^2\backslash\overline{\Omega},
\end{equation}
where $H_0^{(1)}$ is a Hankel function of the first kind of order zero,
and $N$ is the number of approximating functions
each of which is associated with a charge point $\by_j$. It is well
known that $H_0^{(1)}$ satisfies the Helmholtz equation
in $\mathbb{C}\backslash\{0\}$ with a singularity at zero.
It is common to choose charge points lying on a smooth curve;
we then find it enlightening to interpret the MFS
as a discretization of the
single layer potential representation of $u$ as follows.
Let $\Gamma$
be a closed curve enclosing $\overline{\Omega}$ such that
$\text{dist}(\Gamma,\partial\Omega):=\min\{|\bx-\by|,~\bx\in\partial\Omega,
~\by\in\Gamma\}>0$,
then given a density $g\in L^1(\Gamma)$ we may write
\be
u(\bx)\approx\frac{i}{4}\int_\Gamma H_0^{(1)}(k|\bx-\bs|)g(\bs)~\text{d\bs},
\quad\bx\in\Omega.
\label{eq:layer}
\ee
If $g(\bs)=\sum_{j=1}^N\alpha_j\delta(\bs-\by_j)$ for some point set
$\{\by_j\}\in\Gamma$, where
$\delta$ is the Dirac delta, we recover the MFS formulation
\eqref{eq:mfs}.
%
Note that the irregular Bessel function $Y_0$,
or Hankel $H_0^{(2)}$, may be used instead of $H_0^{(1)}$
in the MFS
\cite{En96,Ka01}; see Remark~\ref{rmk:HvsY} below.

An overview about the history of this method and its applications
is given in \cite{FaKa98}. The rate of convergence of the MFS for the Laplace
BVP was investigated in \cite{Bo85,Ka89,Ka90,Ka94,Ka88}. It turns out that
if the boundary data is analytic one can achieve exponential convergence for
the MFS for the Laplace problem on analytic domains if the charge
points $\by_j$ are suitably chosen.

One of the main drawbacks of the method is that in
the end systems of equations or linear least squares problems have to be solved
that are often ill-conditioned. The effects of this ill-conditioning on the
quality of the solution have been investigated for the Laplace problem in
\cite{Ki88,Ki91}.
In this paper we investigate more closely under what conditions on the points
$\by_j$ a numerically stable representation of an approximate solution of
the Helmholtz problem \eqref{eq:bvproblem}
as linear combination of fundamental solutions is possible.
It turns out that this depends on how far into the complex plane a
solution of \eqref{eq:bvproblem} can be analytically continued.
The importance of this in the context of 
scattering problems has already been observed
\cite{Ky85,Ky96} (and references in \cite{Ky96}).

Our work also has consequences for the numerical
solution of more challenging and widely-applicable PDE problems
that are closely related to the one we study. We have in mind
i) finding eigenmodes of the Laplace operator in $\Omega$
with homogeneous boundary conditions
(where the MFS has been used at low \cite{En96} and
very high eigenvalue \cite{que}), and ii)
scattering of time-harmonic waves
(the exterior Helmholtz boundary value problem
in $\mathbb{R}^2\setminus\overline{\Omega}$).
In both these situations
the boundary data is almost always analytic: in problem i) it is homogeneous
and in ii) a plane wave or point source.
%

We will study convergence of the MFS approximation in the boundary error norm
\be
t =
\| \uN - v \|_{L^2(\pO)}\quad .
\label{eq:t}
\ee
By applying \cite[Eq. 7]{Ku78},
this controls the interior error of the solution as follows,
\be
\|\uN - u\|_{L^2(\Omega)}
\le \frac{C_\Omega}{d}\|\uN - v \|_{L^2(\pO)}\quad ,
\label{eq:mp}
\ee
where $d:= \min_j |k^2-E_j|/E_j$,
the domain's Dirichlet eigenvalues are $E_j$,
and $C_\Omega$ is a domain-dependent constant.
This shows that for any fixed nonresonant $k$,
we may use the boundary norm.

In Section \ref{sec:unitdisc} we give rigorous results for the
convergence and the numerical stability of the MFS for Helmholtz
problems on the unit disc, with analytic boundary data,
using charge points on a concentric circle.
We then present a heuristic model for behavior in finite-precision
arithmetic and show it explains well
numerical results observed at both low and high wavenumbers.
A key
conclusion will be that it is the growth in norm
of the coefficient vector
that in practice limits the achievable error,
so this norm should be kept as small as possible to retain high
accuracy.
The reader should take care throughout not to confuse
statements about the {\em coefficient norm}
(which depending on the choice of MFS charge points may either
grow or not grow with $N$ as the error converges to zero),
with statements about the {\em condition number} of the problem
(which always grows with $N$ since the MFS \eqref{eq:mfs}
approximates a single-layer operator \eqref{eq:layer}
which is compact).

In Section \ref{sec:continuation} we move to general analytic domains, and
review results for the analytic continuation of solutions
$u$ of \eqref{eq:bvproblem},
in particular how both the boundary data and the domain
shape may lead to singularities in the continuation of $u$.
In Section \ref{sec:conf} we explore the use of the exterior conformal
map in choosing the charge points $\by_j$
for several more complicated and nonconvex domains.
We propose and provide evidence for conjectures in general domains which are
analogous to the theorems on convergence rate and stability in the unit disc.
In Section \ref{sec:adapt}
we propose and demonstrate a method for choosing charge points
well-adapted to the singularity locations and the wavenumber $k$,
that outperforms the conformal mapping method by a large margin.

\begin{figure} 
\bc
\includegraphics[width=3in]{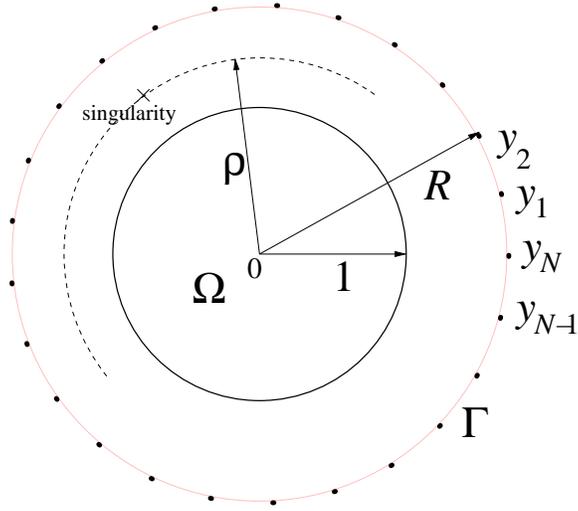}
\ec
\caption{Geometry for the MFS in the unit disc.}
\label{fig:g}
\end{figure} 

\section{The MFS on the unit disc}
\label{sec:unitdisc}

In this section we analyse the accuracy and coefficient
sizes that result in the unit disc $\Omega = \{ \bx: |\bx|< 1\}$,
for the MFS using charge points $\by_j = Re^{i\phi_j}$,
$j=1,2,\ldots,N$, with $\phi_j=2\pi j/N$,
that is, equally spaced on a larger circle of radius $R>1$.
See Figure~\ref{fig:g}. We identify $\mathbb{R}^2$ with $\mathbb{C}$.

Before we embark we need the definition of
the Fourier series for a function $g\in L^2([0,2\pi])$,
\be
g(\theta) = \sum_{m=-\infty}^\infty \hat{g}(m) e^{i m \theta},
\qquad\qquad
\hat{g}(m) = \frac{1}{2\pi}\int_0^{2\pi} g(\theta) e^{-i m \theta} d\theta.
\label{eq:fs}
\ee
Parseval's identity is then
\be
\|g\|_{L^2([0,2\pi])}^2 = 2\pi \sum_{m=-\infty}^\infty |\hat{g}(m)|^2=:2\pi\|\hat{g}\|_{\ell^2(\mathbb{Z})}^2.
\label{eq:pars}
\ee

We also need to represent the coefficient vector
$\mbf{\alpha} := \{\alpha_j\}_{j=1,\ldots, N}$
in a discrete Fourier basis
labeled by $-N/2<k\le N/2$ (we will always choose $N$ even),
\be
\alpha_j = \sum_{k=-N/2+1}^{N/2} \hat{\alpha}_k e^{ik\phi_j},
\qquad\qquad
\hat{\alpha}_k = \frac{1}{N}\sum_{j=1}^{N} \alpha_j e^{-ik\phi_j},
\label{eq:dft}
\ee
where inversion follows from
$\sum_{j=1}^{N}e^{2\pi i kj/N} = N \delta^{(N)}_{k0}$ 
with $\delta^{(N)}_{k0}$, the periodized Kronecker delta defined by
\be
\delta^{(N)}_{kj} = \spl{1,&k\equiv j \pmod N}{0, &\mbox{otherwise}.}
\label{eq:pd}
\ee
Parseval's identity now gives $|\mbf{\alpha}|^2 = N |\hat{\mbf{\alpha}}|^2$,
where $|\mbf{\alpha}|:=\left(|\alpha_1|^2+\dots+|\alpha_n|^2\right)^{1/2}$
is the standard Euclidean norm.

\begin{figure} 
\bc
\mbox{\raisebox{-2.5in}{\includegraphics[height=2.8in]{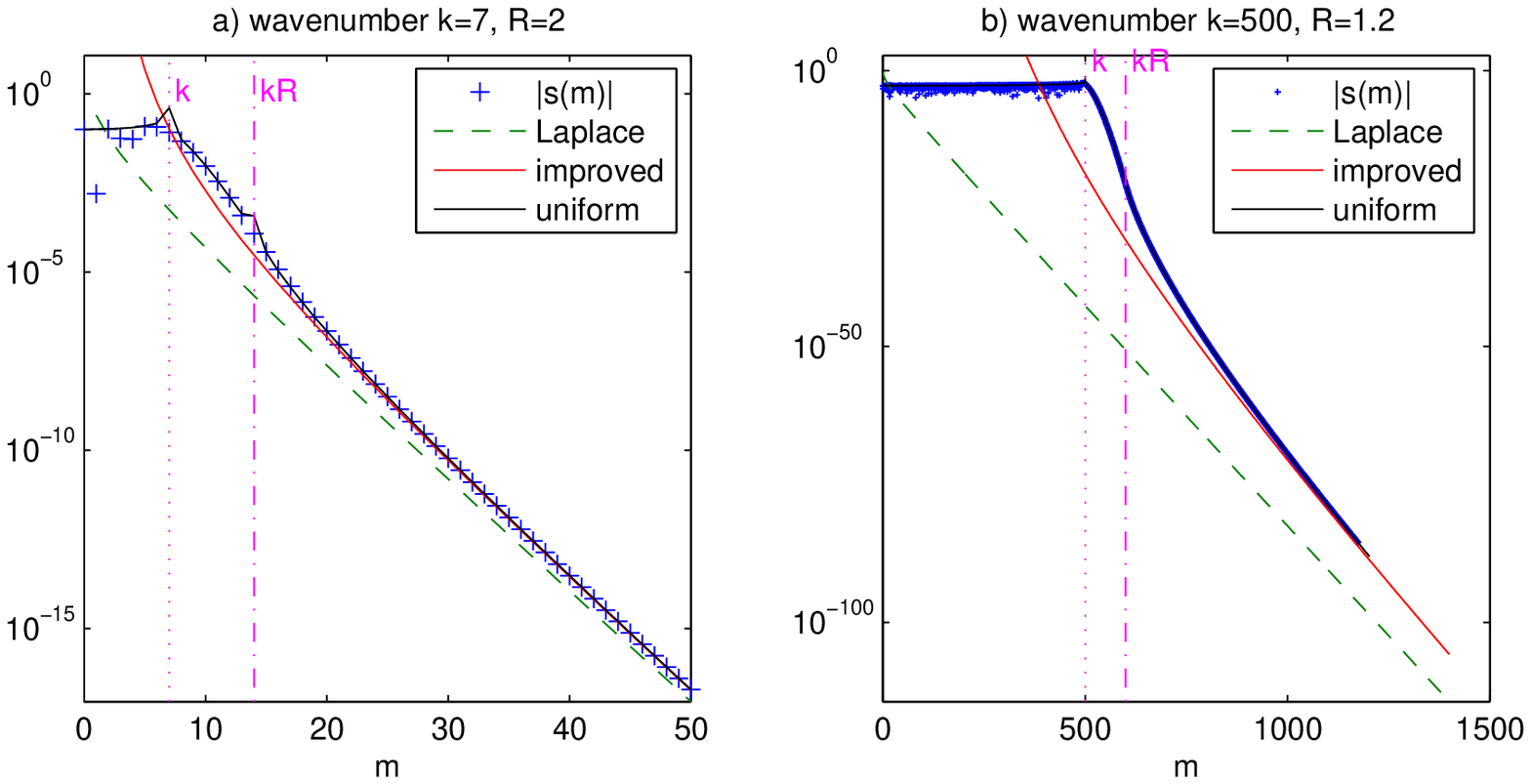}}
c)\raisebox{-2.5in}{\includegraphics[height=2.6in]{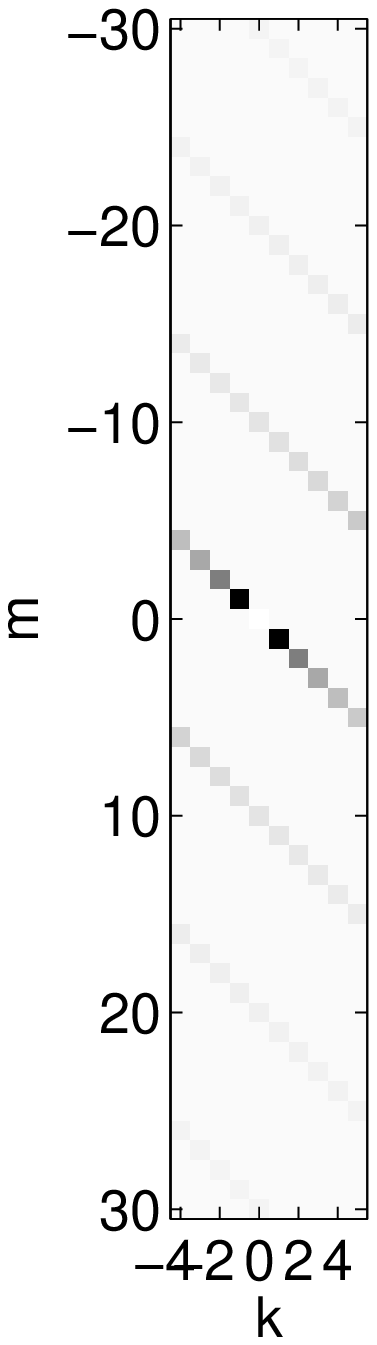}}}
\ec
\caption{Comparison for the unit disc of
layer-potential eigenvalue magnitudes $|\hat{s}(m)|$
given in \eqref{eq:sm} against various
asymptotic expressions: `Laplace' \eqref{eq:lapas},
`improved' \eqref{eq:smser}, and
`uniform' \eqref{eq:unif}.
a) low wavenumber, b) high wavenumber.
Panel c) shows density plot of matrix elements
(\ref{eq:tmk}) of $Q$ for $N=10$,
in the domain $|m|\le 30$. In c) we
chose unrealistically small values of $R$ and $k$
in order to make the super- and sub-diagonals more visible.
}
\label{fig:sm}
\end{figure} 

\subsection{Map from layer potential to Fourier basis on the unit circle}

For simplicity we first consider the layer potential version
of this problem, which can be interpreted as the
$N\to\infty$ limit of the MFS.
The single layer potential lying on the outer circle
$\Gamma = \{\by : |\by|=R\}$ is 
\be
u(\bx) = 
\frac{i}{4}\int_0^{2\pi} H_0^{(1)}(k|\bx-Re^{i\phi}|)g(\phi) \, d\phi.
\label{eq:uinf}
\ee
Note $g\in L^1([0,2\pi])$ is the density with respect to angle measure
$d\phi$ rather than the usual length measure $Rd\phi$.
The Fourier-Bessel decomposition of a fundamental solution located
at $Re^{i\phi}$, evaluated at $\bx=re^{i\theta}$ is,
using Graf's addition formula \cite[Eq. 9.1.79]{a+s},
\be
\frac{i}{4} H_0^{(1)}(k|\bx - Re^{i\phi}|) =
\frac{i}{4} \sum_{m\in\mathbb{Z}} H_m^{(1)}(kR) J_m(kr) \cos m(\theta-\phi)
=
\frac{i}{4} \sum_{m\in\mathbb{Z}} H_m^{(1)}(kR)e^{-im\phi} \cdot
J_m(kr) e^{im\theta}
\label{eq:fund}
\ee
where the second step involved the reflection formulae \cite[Eq. 9.1.5]{a+s}
$J_{-m}(z)=(-1)^m J_m(z)$ and $Y_{-m}(z) = (-1)^m Y_m(z)$.
Hence the Fourier-Bessel coefficients are
$\frac{i}{4}H_m^{(1)}(kR)e^{-im\phi}$.
The restriction of (\ref{eq:uinf}) to $\bx\in\pO$ gives a single-layer
operator $S:L^2([0,2\pi])\to L^2([0,2\pi])$
of convolution type, which is therefore
diagonal in the Fourier basis $\{e^{im\theta}\}_{m\in\mathbb{Z}}$
and entirely described by its eigenvalues.
Comparing (\ref{eq:uinf}), (\ref{eq:fund})  and using orthogonality
gives
$\hat{u}(m)  = \hat{s}(m) \hat{g}(m)$ where
the eigenvalues of $S$ are
\be
\hat{s}(m) = \frac{i\pi}{2} H_m^{(1)}(kR) J_m(k).
\label{eq:sm}
\ee
\begin{rmk}        

Since the Hankel function (real argument) is never zero,
an eigenvalue can vanish only when
$J_m(k)=0$, corresponding to
a Dirichlet eigenvalue (resonance) of $\Omega$.
In contrast if $Y_0$ were chosen as the fundamental solution
in \eqref{eq:mfs}, $Y_m(kR)$ may accidentally be very small giving
poor or spurious numerical
results, although in practice this happens rarely \cite{En96}.
In general one may avoid this problem by using
$(i/4)(Y_0 + i\eta J_0)$, for any real $\eta\neq 0$
(an analogous idea is used in layer potentials, p.48 of \cite{coltonkress}).
\label{rmk:HvsY}
\end{rmk}            
Since its kernel
is continuous $S$ is compact, so
$\lim_{|m|\to\infty}\hat{s}(m) = 0$.
The compactness of $S$ means its inverse is unbounded, and
we expect to find arbitrarily large $\|g\|$ needed to
represent certain unit-norm boundary functions $v$.

In the Laplace ($k\to0$) limit we recover the following known result
(\eg \cite[Eq 3.2]{Ka90}).
We use the small-argument asymptotics
$J_m(k)\sim (k/2)^m/\Gamma(m+1)$
and $Y_m(kR)\sim-\frac{1}{\pi}\Gamma(m)
(kR/2)^{-m}$ for all integer $m>0$, and the reflection
formulae, and get
\be
\hat{s}(m) \rightarrow
\frac{1}{2|m|} R^{-|m|}, \qquad m \in \mathbb{Z}\backslash\{0\},
\quad k\rightarrow 0. 
\label{eq:smlap}
\ee
To analyze convergence rate we will need the
asymptotic behavior as $|m|\to\infty$ for fixed $k$.
Using in \eqref{eq:sm} the large-order asymptotics
(9.3.1 in \cite{a+s})
$J_m(z)\sim \frac{1}{\sqrt{2\pi m}}(ez/2m)^m$
and
$Y_m(z)\sim -\sqrt{\frac{2}{\pi m}}(ez/2m)^{-m}$, where $z$ is fixed,
gives the leading-order behavior
\be
\hat{s}(m) \sim \frac{1}{2|m|} R^{-|m|}, \qquad |m|\to\infty,
\label{eq:lapas}
\ee
which coincides with the Laplace case \eqref{eq:smlap}.
We therefore have the exponential uniform bounds,
for some constants $c_s$ and $C_s$ depending only on $R$ and $k$,
\be
\frac{c_s}{|m|}
R^{-|m|} \le |\hat{s}(m)| \le \frac{C_s}{|m|} R^{-|m|}\leq
C_s R^{-m}, \qquad m\in\mathbb{Z}\backslash\{0\},
\label{eq:smbnd}
\ee
where $C_s$ is chosen large enough such that also $|\hat{s}(0)|\leq C_s$.
Figure~\ref{fig:sm}a shows that for low wavenumbers the leading-order
asymptotic is reached rapidly, hence
the smallest possible ratio $C_s/c_s$
is not too large (here about $10^3$ will suffice).
This asymptotic becomes accurate well before the dynamic
range begins to exceed machine precision ($\emach\approx 10^{-16}$
for double precision).
However, the situation can differ radically for large wavenumbers,
as Figure~\ref{fig:sm}b illustrates. Here the Laplace asymptotic
is not relevant for the eigenvalues within a factor $\emach$ of the
maximum. Worse still, the ratio $C_s/c_s$
must be exceedingly large (several
tens of orders of magnitude).
We will present more useful
asymptotic approximations for the eigenvalues
in Section \ref{sec:est}.
For now we need only \eqref{eq:smbnd} to prove
exponential convergence rates.

\subsection{Map from MFS coefficients to Fourier basis}

We now adapt the above to the discrete source case. Define the density
\be
g
(\phi) = \sum_{j=1}^{N} \alpha_j \delta(\phi-\phi_j).
\label{eq:gmfs}
\ee
It follows that
$$
u^{(N)}(e^{i\phi})=\frac{i}{4}\sum_{j=1}^N
\alpha_jH_0^{(1)}(k|e^{i\phi}-Re^{i\phi_j}|)
=(Sg
)(\phi).
$$
We have
$$
\|Sg
\|_{L^2([0,2\pi])}^2=2\pi\|\widehat{Sg
}\|_{\ell^2(\mathbb{Z})}^2=2\pi\sum_{m=-\infty}^{\infty}
|\hat{s}(m)\hat{g}
(m)|^2
=\frac{N^2}{2\pi}
\sum_{m=-\infty}^{\infty}|\hat{s}(m)\hat{\alpha}_{m\bmod N}|^2,
$$
where $m\bmod N$ denotes the unique integer lying in the range
$-N/2+1, \ldots, N/2$ which differs from $m$ by an integer multiple of $N$.
The last equality follows from the Fourier series representation of
\eqref{eq:gmfs},
\begin{equation}
\hat{g}
(m) = \frac{1}{2\pi}\sum_{j=1}^{N} \alpha_j
e^{-im\phi_j}= \frac{N}{2\pi} \hat{\alpha}_{m\bmod N },
\qquad m\in\mathbb{Z}.
\label{eq:gfourier}
\end{equation}

Applying H\"older's inequality
and \eqref{eq:smbnd} we obtain
\be
\|Sg
\|_{L^2([0,2\pi])}^2
\leq
\frac{N^2}{2\pi}\max_{j=-\frac{N}{2}+1\dots\frac{N}{2}}
|\hat{\alpha}_j|^2
\cdot\!\!
\sum_{m=-\infty}^{\infty}|\hat{s}(m)|^2\leq
\frac{N^2}{2\pi}|\hat{\mbf{\alpha}}|^2 C_s^2\frac{R^2+1}{R^2-1}.
\label{eq:sgnorm}
\ee

Define the operator $Q:\mathbb{R}^N \to \ell^2(\mathbb{Z})$ by
$Q\hat{\mbf{\alpha}}:=\widehat{Sg
}$.
Therefore $Q$ maps the discrete Fourier coefficient vector
$\hat{\mbf{\alpha}}$ to the Fourier series coefficients on the boundary $\pO$.
From \eqref{eq:sgnorm} we immediately obtain the following.

\begin{lemma} 
For $R>1$ the operator $Q$ is bounded. Furthermore,
$$
\|Q\|
\leq
C_s\frac{N}{2\pi}\sqrt{\frac{R^2+1}{R^2-1}},
$$
where 
$
\|Q\|:=\max_{\hat{\mbf{\alpha}}\in\mathbb{R}^N\backslash\{0\}}\frac{\|Q\hat{\mbf{\alpha}}\|_{\ell^2(\mathbb{Z})}}{|\hat{\mbf{\alpha}}|}$.
\end{lemma} 

The action of $Q$ is that of a generalized
matrix of width $N$ but (bi-)infinite height,
\be
\hat{u}(m) = \sum_{k=-N/2+1}^{N/2} q_{mk} \hat{\alpha}_k,
\qquad \mbox{for} \; m\in\mathbb{Z}.
\label{eq:um}
\ee
From $Q\hat{\mbf{\alpha}}(m)=(\widehat{Sg
})(m) =
\hat{s}(m)\hat{g}
(m)$ and \eqref{eq:gfourier} it follows that
the matrix elements are
\be
q_{mk} = \frac{N}{2\pi} \hat{s}(m) \delta^{(N)}_{mk} .
\label{eq:tmk}
\ee
Fig.~\ref{fig:sm}c) shows a greyscale picture of a piece
of the resulting matrix $Q$.
Notice that it is dominated by a main diagonal proportional to the
diagonal of the $S$ operator defined in (\ref{eq:sm}), but with
(exponentially) smaller entries
on an infinite sequence of super- and sub-diagonals.
This off-diagonal part can be interpreted as aliasing `overtones'
due to discrete sampling of a continuous layer potential.
In the Laplace case using \eqref{eq:smlap} in \eqref{eq:tmk}
recovers the results of Katsurada \cite[Lemma 1, case 2]{Ka88}.

\subsection{Convergence rate and coefficient sizes
in the disc with analytic data}

We are now in a position to express the boundary error norm \eqref{eq:t}
in terms of $\hat{\mbf{\alpha}}$.
Combining with \eqref{eq:pars} and \eqref{eq:um} gives
\be
t[\hat{\mbf{\alpha}}] = \sqrt{2\pi}\;
\| Q \hat{\mbf{\alpha}} - \hat{\mbf{v}}\|_{\ell^2(\mathbb{Z})}\quad ,
\label{eq:tlsq}
\ee
where $\hat{\mbf{v}}\in \ell^2(\mathbb{Z})$ is the
p a recipe for locating charge points which
allows us to reach error norms of typically $10^{-11}$ on a wide variety
of analytic domains.
At high frequencies of order only 3 points per wavelength are
needed, which compares very favorably to boundary integral methods.

vector of Fourier coefficients of the boundary data $v$ on the unit circle.
Assume that $v$ can be analytically continued to
the annulus $\{z\in\mathbb{C}:~\frac{1}{\rho}<|z|<\rho\}$ for some $\rho>1$, that is the closest singularity of the analytic continuation of $\rho$ has the radius $\rho$ or $1/\rho$. We then have asymptotically exponential decay of the Fourier coefficients,
\be
%
|\hat{v}(m)| \;\sim\; C \rho^{-|m|}, \qquad |m|\to\infty,
\label{eq:decay}
\ee
for some constant $C$.
A simple example is boundary data arising
from an $n^{th}$-order pole $v(z) = \mbox{Re}\,(z-\rho)^{-n}$ for
$z\in\pO$, $n=1,2,\ldots$~.

Minimizing \eqref{eq:tlsq}
over $\hat{\mbf{\alpha}}$ is a least-squares problem
involving the generalized matrix $Q$.
But since the columns of $Q$ are orthogonal this separates into $N$
independent single-variable minimizations.
We may use a diagonal approximation to choose $\hat{\mbf{\alpha}}$
which is sufficient for the following convergence rate bounds.
\begin{thm} 
Let  $R>1$ and $N$ be even. 
For analytic boundary data $v$ obeying (\ref{eq:decay}), the
minimum boundary error \eqref{eq:t} achievable with the MFS
in the unit disc satisfies
\be
t \; \le \; \left\{\begin{array}{ll}
C\rho^{-N/2},& \rho < R^2\\
C\sqrt{N}R^{-N},& \rho = R^2\\
CR^{-N},& \rho > R^2
\end{array}\right.
\ee
where each time $C$ means a different constant
which may depend on $k$, $R$, and $v$, but not $N$.
Furthermore if $v$ is analytically continuable to an entire function,
the last of the three cases holds for any $R>1$.
\label{thm:t}
\end{thm} 
\bp
We choose coefficients
$\hat{\alpha}_m = \hat{v}(m)/q_{mm}$ for $-N/2<m\le N/2$.
This exactly matches the Fourier coefficients in this interval,
therefore errors are due only to frequencies lying outside the interval.
\eqref{eq:tlsq}, \eqref{eq:tmk} and the triangle
inequality in $\ell^2(\mathbb{Z})$ give
\be
t \; = \; \left(2\pi\sum_{m\notin[-\frac{N}{2}+1,\frac{N}{2}]} \left|
(Q\hat{\mbf{\alpha}})(m) - \hat{v}(m)\right|^2 \right)^{1/2}
\; \le \; \sqrt{2\pi} (E_u + E_v),
\nonumber
\ee
where
\be
E_u^2 \;=\; \sum_{m\notin[-\frac{N}{2}+1,\frac{N}{2}]}
|(Q\hat{\mbf{\alpha}})(m)|^2
\;=\;
\sum_{-\frac{N}{2}<n\le \frac{N}{2}}
\left| \frac{\hat{v}(n)}{\hat{s}(n)}\right|^2
\sum_{b\neq 0}\left| \hat{s}(bN+n)\right|^2
\label{eq:eusum}
\ee
and
\be
E_v^2 \; = \; \sum_{m\notin[-\frac{N}{2}+1,\frac{N}{2}]}
|\hat{v}(m)|^2\quad .
\ee
We can bound both error terms since all terms in the sums
have exponential bounds.
First we note that
using \eqref{eq:smbnd} and \eqref{eq:decay} gives
\begin{align}
E_u^2&\leq C_1\sum_{\substack{-\frac{N}{2}<n\le \frac{N}{2}\\ n\neq 0}}
\left(\frac{R}{\rho}\right)^{2|n|}|n|^2\sum_{b\neq 0} \frac{R^{-2|bN+n|}}{|bN+n|^2}+
C_2\left|\frac{\hat{v}(0)}{\hat{s}(0)}\right|^2\sum_{b\neq 0} \frac{R^{-2|bN|}}{|bN|^2}\nonumber\\
&\leq C_1\sum_{\substack{-\frac{N}{2}<n\le \frac{N}{2}\\ n\neq 0}}
\left(\frac{R}{\rho}\right)^{2|n|}\sum_{b\neq 0} R^{-2|bN+n|}+
C_2\left|\frac{\hat{v}(0)}{\hat{s}(0)}\right|^2\sum_{b\neq 0} R^{-2|bN|}
\label{eq:smsum}
\end{align}
for sufficiently large constants $C_1$ and $C_2$.
We can bound
\be
\sum_{b\neq 0} R^{-2|bN+n|} \;\le\; C_3 R^{-2N+2|n|},\quad
-\frac{N}{2}<n\leq \frac{N}{2}
\label{eq:bsum}
\ee
for a large enough constant $C_3$. Inserting \eqref{eq:bsum} into
\eqref{eq:smsum} and absorbing
$\left|\frac{\hat{v}(0)}{\hat{s}(0)}\right|^2$ into the constants gives
\be
E_u^2   \;\le\;
CR^{-2N}\sum_{-\frac{N}{2}<n\le \frac{N}{2}}
\left(\frac{R^2}{\rho}\right)^{2|n|}
\label{eq:eubnd}
\ee
for a sufficiently large constant $C$.
Similarly $E_v^2\le C \rho^{-N}$ follows from \eqref{eq:decay}.
We now study the sum in (\ref{eq:eubnd}).
For $\rho<R^2$, \eqref{eq:eubnd} can be estimated by $C\rho^{-N}$ for
some constant $C>0$.
This means both
$E_u^2$ and $E_v^2$ have the same exponential decay $C\rho^{-N}$.
For $\rho=R^2$, the sum contains $N$ equal terms and therefore
$E_u^2\leq CNR^{-2N}$, which decays slower than the bound $C\rho^{-N}$
on $E_v^2$.
For $\rho>R^2$, \eqref{eq:eubnd} can be estimated by $CR^{-2N}$ for
some $C>0$, which means $E_v$
is of higher negative order in $N$ than $E_u$ and can be dropped.
For the case of $v$ continuable to an entire function,
we may take $\rho\to\infty$ and the case $\rho>R^2$ applies.
\ep

This is a generalization of a result of Katsurada
\cite{Ka89} from the Laplace to the Helmholtz problem.
Since only the exponential bounds \eqref{eq:smbnd} and no other
information about $\hat{s}(m)$ was used, the convergence
rates are identical to those for Laplace with the same boundary data.


\begin{rmk}
An interpretation of the two main
convergence rate regimes is,
\bi
\item $v$ is `not relatively smooth' ($R>\sqrt{\rho}$, \ie `distant' charge points):
errors are limited by the absence
of Fourier modes beyond a frequency $N/2$ in the MFS basis,
hence rate is controlled by the boundary data singularity $\rho$.
\item $v$ is `relatively smooth' ($R<\sqrt{\rho}$, \ie `close' charge points):
errors are limited by
aliasing errors due to the discrete representation of the single
layer potential, hence rate is controlled by $R$.
\ei
\label{rmk:interp}
\end{rmk}

\begin{rmk}
By keeping track of the constants in the above proof, one can check that
$C$ is at least as big as $C_s/c_s$, which, as discussed, must be
very large for large $k$.
Thus while the convergence rates in the theorem must be reached
asymptotically as $N\to\infty$ in exact arithmetic, we cannot expect the
bounds to be numerically useful in practice at large wavenumbers.
\end{rmk}

We are interested in how the coefficient norm $|\mbf{\alpha}|$
grows as we reduce the boundary error in the MFS representation
\eqref{eq:mfs} by increasing $N$.
Firstly, it is easy to show that when the MFS charge points are
closer than the nearest singularity, the coefficients need not grow.
\begin{thm} 
Let $\Omega$ be the unit disc, and $R<\rho$,
with fixed analytic boundary data $v$ obeying (\ref{eq:decay}).
Then as $N\to\infty$
there exists a sequence of coefficient vectors $\bmal$
with bounded norm $|\bmal|$,
with corresponding boundary error norm \eqref{eq:tlsq} converging as in
Thm.~\ref{thm:t}.
\label{thm:n}
\end{thm} 
\bp
We choose coefficients as in the proof of Thm.~\ref{thm:t}, which
therefore give the desired convergence rate.
Then using \eqref{eq:decay}, \eqref{eq:tmk} and \eqref{eq:lapas},
\be
\hat{\alpha}_m = \frac{\hat{v}(m)}{q_{mm}} =
\frac{2\pi}{N}\frac{\hat{v}(m)}{\hat{s}(m)}
\sim
C \frac{|m|}{N} \left(\frac{R}{\rho}\right)^{|m|}\!\!\!
\leq
\frac{C}{2}\left(\frac{R}{\rho}\right)^{|m|}\!\!\!,
\quad
-\frac{N}{2}<m\le\frac{N}{2},~m\neq 0
\label{eq:diag}
\ee
for some constant $C$.
For $R<\rho$ this is an exponentially decaying sequence so, independent of
$N$, $|\bmal|$ is bounded by a constant.
\ep

More problematically, the coefficient choice used in the above two proofs
would then imply in the case $\rho<R$
that
$|\mbf{\alpha}|$ diverges exponentially with $N$.
However, it is not immediately obvious whether there is a
different
choice of $\hat{\alpha}_m$ which avoids exponential growth.
The following theorem excludes this possibility by showing that
when the singularity in the analytic continuation of the boundary
data is closer than the MFS source points,
{\em any} convergent sequence of coefficient vectors $\mbf{\alpha}$
must diverge in norm in this way.
\begin{thm} 
Let $\Omega$ be the unit disc, with $R>\rho$.
Let the boundary data Fourier coefficients decay no faster than
\eqref{eq:decay}, that is, for some constant $c_v$,
\be
|\hat{v}(m)| \ge c_v \rho^{-|m|}.
\label{eq:vlb}
\ee
For any positive even $N$ satisfying
$N > 3+N_\tbox{min}$, where
$N_\tbox{min}:=2\max\left[\ln{\left(\sqrt{\frac{2}{\pi}}
\frac{c}{c_v}\right)} / \ln \rho,\,1\,\right]$,
let $\mbf{\alpha}$ be a coefficient vector
such that the MFS representation \eqref{eq:mfs} has a
boundary error norm \eqref{eq:tlsq} satisfying
\be
t \le c \rho^{-N/2},
\label{eq:tconv}
\ee
where $c$ is a constant independent of $N$.
Then
\be
|\mbf{\alpha}| \ge C\sqrt{N}\left(\frac{R}{\rho}\right)^{N/2}
\ee
for some constant $C$ which may depend on $k$, $R$, and $v$, but not $N$.
\label{thm:a}
\end{thm} 
Note that \eqref{eq:tconv} is the appropriate convergence rate
for the case $\rho<R^2$ derived in Theorem~\ref{thm:t}.

\bp 
For even $N$ fix $N > N_\tbox{min}+3$.
Using \eqref{eq:tmk} in \eqref{eq:tlsq} implies the trivial bound
\be
\sqrt{2\pi}
\left| \frac{N}{2\pi}\hat{s}(m)\hat{\alpha}_{m \bmod N} - \hat{v}(m)\right|
\le t, \qquad \mbox{for all } m\in\mathbb{Z}.
\label{eq:triv}
\ee
Define the (positive) maximum Fourier frequency $F:=
\frac{N}{2}-K$, where $K$ is the unique integer such that
$$
\frac{N_\tbox{min}}{2}\leq K<\frac{N_\tbox{min}}{2}+1.
$$
Note that $K$ is independent of $N$.
One can verify using
\eqref{eq:tconv}, \eqref{eq:vlb}
and the definition of $N_\tbox{min}$ that
\be
|\hat{v}(m)| \ge \sqrt{\frac{2}{\pi}}t,
\qquad \mbox{ for all } |m|\le F.
\label{eq:F}
\ee
In the frequency range $|m|\le F$, it follows from \eqref{eq:triv}
and \eqref{eq:F} that $\hat{v}(m)$ is sufficiently large relative to
$t$ to bound the coefficients away from zero,
\be
|\hat{\alpha}_m| \;\ge\; \frac{2\pi}{N}
\frac{|\hat{v}(m)|- t/\sqrt{2\pi}}{|\hat{s}(m)|}
\;\ge\; \frac{\pi}{N}\frac{|\hat{v}(m)|}{|\hat{s}(m)|}
\;\ge\; \frac{\pi c_v|m|}{C_s N}\left(\frac{R}{\rho}\right)^{|m|},\qquad
0<|m|\leq F
\ee
where the last step used \eqref{eq:smbnd} and \eqref{eq:vlb}.
Choosing the maximal frequency $m=F=\frac{N}{2}-K$ we obtain
\be
|\hat{\alpha}_F|\;\geq\;
\frac{\pi
  c_v}{C_s}\left(\frac{R}{\rho}\right)^{\frac{N}{2}-K}\!
\left(\frac{1}{2}-\frac{K}{N}\right)
\;>\;
\frac{\pi
  c_v}{2C_s}\left(\frac{R}{\rho}\right)^{\frac{N}{2}-K}
\!\!\!\frac{1}{N_\tbox{min}+3}.
\label{eq:uglyest}
\ee
Here the latter inequality follows from
$$
\frac{1}{2}-\frac{K}{N}\;>\;
\frac{1}{2}-\frac{N_\tbox{min}/2+1}{N_\tbox{min}+3}\;=\;
\frac{1}{2(N_\tbox{min}+3)}.
$$

Absorbing the $N$-independent factors of \eqref{eq:uglyest} into a
constant and
noticing that the Euclidean norm of a vector is at least as large as
its largest component we have
\be
|\mbf{\alpha}| = \sqrt{N} |\hat{\mbf{\alpha}}| \;\ge \;
\sqrt{N}|\hat{\alpha}_F|\;\ge\;
C\sqrt{N}\left(\frac{R}{\rho}\right)^{N/2}   
\ee
for a sufficiently small constant $C>0$.
\ep 

\begin{figure} 
\center
\includegraphics[width=0.8\textwidth]{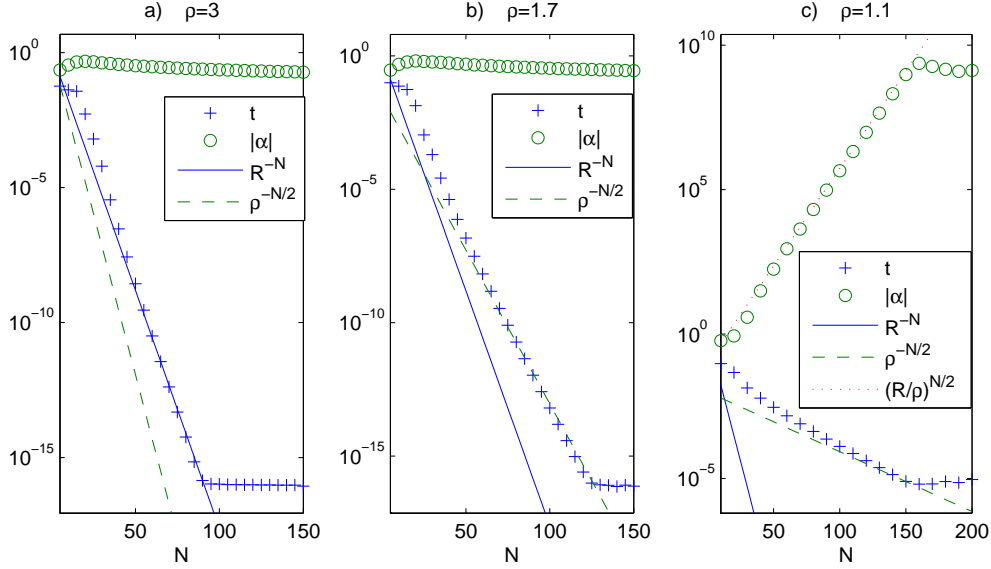}
\caption{Convergence and coefficient sizes as a function of $N$, for
the MFS approximation to the interior Helmholtz BVP in the unit disc
given boundary data corresponding to a single source
\eqref{eq:vfund} at radius $\rho$.
The wavenumber is low ($k=8$).
The MFS sources are at $R=1.5$.
For visual comparison the relevant power laws
from Theorems~\ref{thm:t} and \ref{thm:a}
are shown
(sometimes the constants have been chosen to match the data).
There were $M=240$ boundary points.
}\label{fig:lowk}
\end{figure} 

\begin{figure} 
\center
\mbox{a)\raisebox{-2.3in}{\includegraphics[width=2.5in]{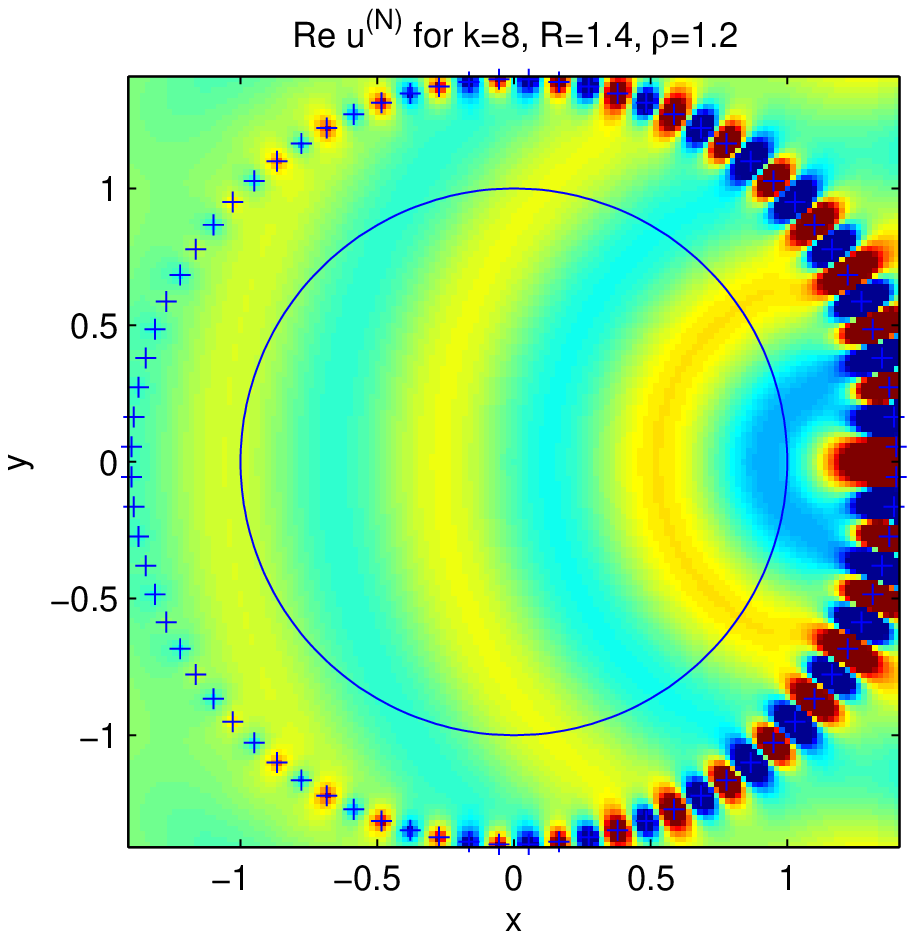}}
\qquad
b)\raisebox{-2.3in}{\includegraphics[width=2.5in]{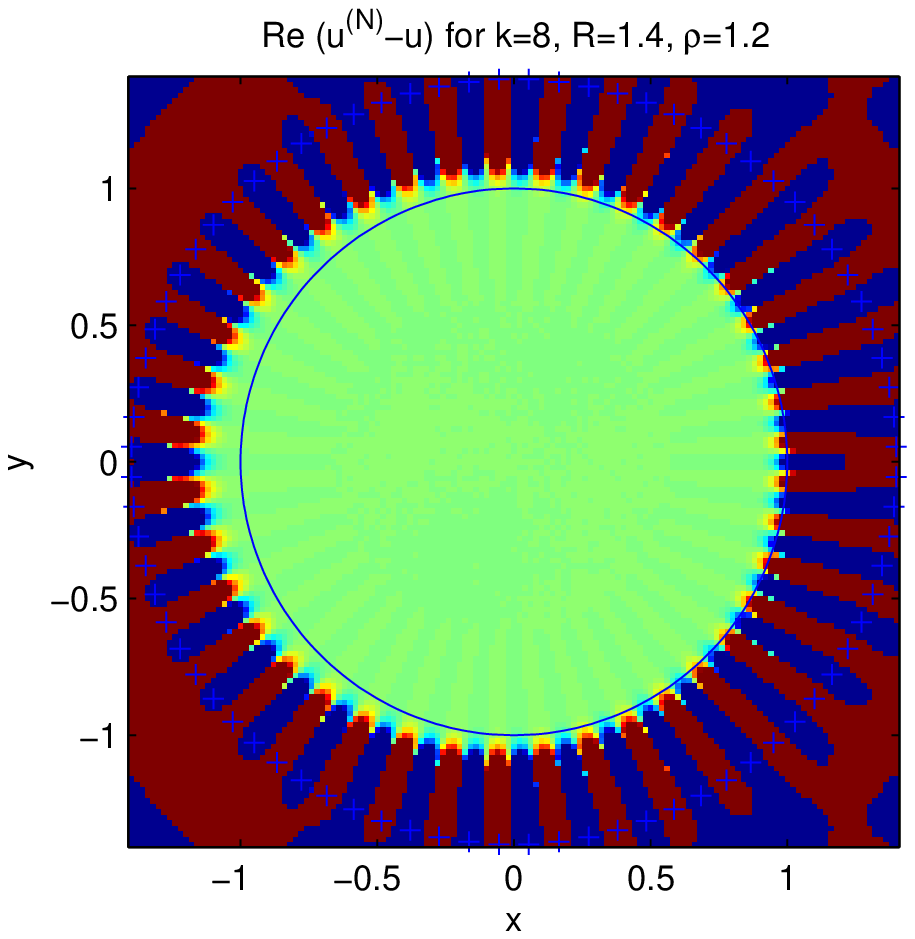}}
}
\caption{MFS approximation to the interior Helmholtz BVP
in the disc with given boundary data corresponding to a single source
\eqref{eq:vfund} outside the domain at $\rho=1.2$, with $N=80$ and $R=1.4$.
a) Re $\uN$, b) residual Re $(\uN - u)$ on a color scale $3\times10^4$
times more sensitive than a).
}\label{fig:udisc}
\end{figure} 

An immediate consequence is that if $\rho$ is known for given
boundary data,
to prevent exponential growth in coefficients one should
restrict the charge point radius to $R\le \rho$.
%
To illustrate this and Thm~\ref{thm:t}, we
finish this section with some numerical experiments at low wavenumber.
The implementation was standard, as follows.
The integral in \eqref{eq:t} is approximated using uniform quadrature with
$M$ equally-spaced boundary points $\{\bx_m\}_{m=1\cdots M}$.
Specifically the $M$-by-$N$ matrix $A$ has elements
\be
A_{mj} : = \frac{i}{4} H^{(1)}_0 (k|\bx_m - \by_j|),
\label{eq:matrix}
\ee
and the boundary-value vector $\mbf{v}\in\mathbb{C}^M$ has elements
$v_m := v(\bx_m)$.
The resulting linear system (usually overdetermined, $M>N$, in our work)
$A\bmal = \mbf{v}$ was solved in the least-squares sense via
the QR decomposition
(MATLAB's backslash command) in double-precision arithmetic.
The boundary error norm then is approximately
$t = \sqrt{|\pO|/M}\, |A\bmal - \mbf{v}|$.

In Fig.~\ref{fig:lowk} we show convergence of $t$ using boundary data
\be
v(z) = -\frac{1}{4} Y_0(k|z-\rho|),   \qquad z\in\pO
\label{eq:vfund}
\ee
with real $\rho>1$, that is, a single real-valued fundamental solution.
The three panels illustrate the three cases of Theorem~\ref{thm:t}.
In a) $\rho>R^2$ thus convergence is determined by $R$.
In b) $R^2>\rho>R$ so we have transitioned to a convergence rate
given by $\rho$.
In both these cases the coefficient size
$|\mbf{\alpha}|$ is very close to constant
(note by contrast that the condition number of $A$ is growing).
However in c) $\rho<R$ so convergence rate is again determined
by $\rho$, but now $|\mbf{\alpha}|$ grows exponentially
at precisely the rate indicated by Theorem~\ref{thm:a}.
Fig.~\ref{fig:udisc} shows the resulting approximate
field $u^{(N)}$ in the case c), and the error function $u^{(N)}-u$.
Notice that the error is oscillatory at Fourier frequencies
of about $N/2$ (as predicted by Remark \ref{rmk:interp};
this can be seen by comparing
the alternating signs in Fig.~\ref{fig:udisc}b) to the angular spacing
of source points),
is concentrated on the side of $\pO$ nearest the singularity,
and decays exponentially inside the domain (it is evanescent).
We have also substituted $v(z) = \mbox{Re}(z-\rho)^{-1}$ and find the
convergence rates in Fig.~\ref{fig:lowk} are very similar.
We note that in each plot in this figure, the convergence eventually
stops, as we now explain.

\begin{figure} 
\center
\includegraphics[width=0.8\textwidth]{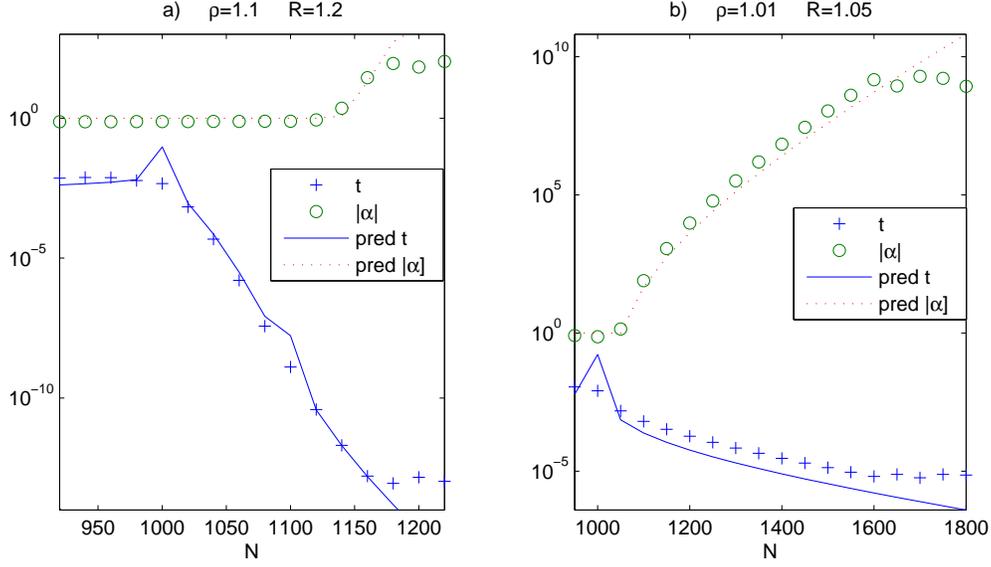}
\caption{Convergence and coefficient sizes as a function of $N$, for
MFS approximation to the interior Helmholtz BVP in the disc
at high wavenumber $k=500$.
The boundary data corresponds to a single source
\eqref{eq:vfund} outside the unit disc at radius $\rho$,
with MFS sources at $R$. For a) $\rho=1.1$, $R=1.2$, b) $\rho=1.01$,
$R=1.05$.
The `predicted' curves are given by \eqref{eq:tap} using \eqref{eq:nap},
with \eqref{eq:unif} modeling both $\hat{s}(m)$ and $\hat{v}(m)$.
}\label{fig:highk}
\end{figure} 

\begin{figure} 
\center
\includegraphics[width=0.8\textwidth]{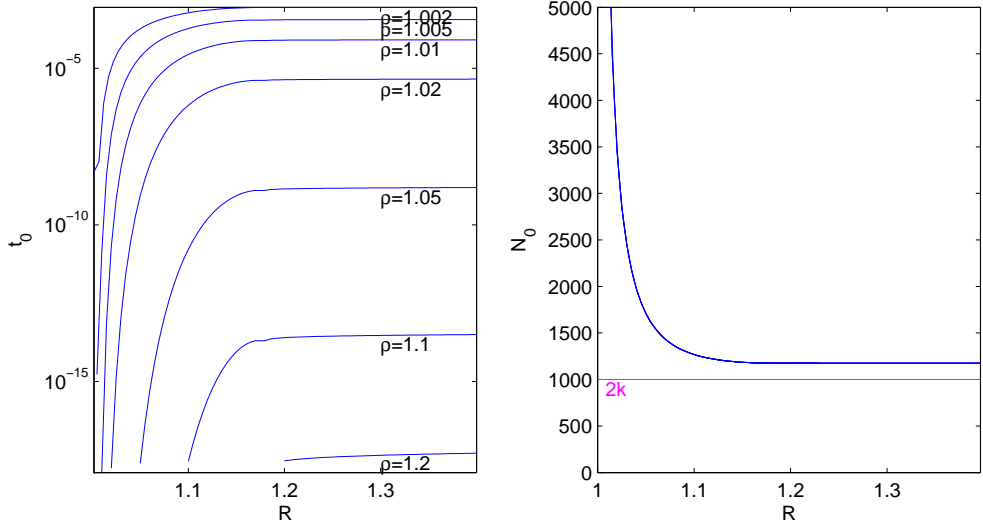}
\caption{a) (left plot) Minimum achievable boundary error norms $t_0$
and b) (right plot) corresponding basis sizes $N_0$, predicted for a selection
of source radii $\rho$ and
MFS source point radius $R$, for the Helmholtz BVP in the unit disc
at high wavenumber $k=500$.
The predictions are done in the range $R>\rho$ using the model
in Sec.~\ref{sec:est}.
Note that in b) all the graphs for different $\rho$ lie on top of one another.
}\label{fig:Rsweep}
\end{figure} 

\subsection{Minimum achievable error in the disc for low and high wavenumbers}
\label{sec:est}

So far we have proven results which hold in exact arithmetic.
However machine precision limits the dynamic range of
eigenvalues that may be used: since
the MFS trial functions in \eqref{eq:mfs}
have typical size of $O(1)$ in $\Omega$ (for any reasonable wavenumber),
each coefficient $\alpha_j$ will
result in round-off errors
of size roughly $\emach \alpha_j$ in the numerical approximation
$u^{(N)}$.
Thus we expect convergence to stop when $t$ reaches of order
$\emach$ times the coefficient norm $|\mbf{\alpha}|$.
This behavior is well illustrated in the three plots of
Fig.~\ref{fig:lowk}: convergence stops when the ratio between $t$
and $|\mbf{\alpha}|$ reaches roughly $10^{-16}$.
In c) the coefficient growth
thus limits achievable error norm to only about $10^{-5}$.
Such premature
halting of convergence has been observed in the Laplace ($k=0$) case
in the disc \cite{Ka89,Ka96,Smyrl} but not analyzed much before.
We analyse this in the Helmholtz case after observing the following
consequence of Thm.~\ref{thm:t} and \ref{thm:n}.

\begin{rmk} 
For any wavenumber, for boundary data with a given singularity radius
$\rho$, the choice of charge point radius $R$ in the range
$\sqrt{\rho}<R<\rho$ leads to both
optimal asymptotic convergence rate $t\sim\rho^{-N/2}$ and a lack of
coefficient growth.
\label{rmk:bestR}
\end{rmk} 

It is useful to have a heuristic model which predicts,
for general $\rho$ and $R$ in the unit disc,
both the
lowest achievable error norm and the basis size $N$ required to achieve it.
We spend the rest of this section constructing then testing such a model.
We first consider low wavenumbers, that is, ones where the Laplace
asymptotic form \eqref{eq:smlap} is relevant for the relevant eigenvalues
(those no smaller than $\emach$ times the largest eigenvalue).
We discuss both the case $R<\rho$ where no coefficient growth occurs,
and the case $R>\rho$ when the coefficient norm grows. 
The proof of Theorem~\ref{thm:a} suggests that, at least
when $\hat{s}(m)$ is exponentially decaying,
the diagonal approximation for the MFS Fourier coefficients
(see \eqref{eq:diag})
\be
\hat{\alpha}_m \approx \frac{2\pi}{N} \,\frac{\hat{v}(m)}{\hat{s}(m)}
\ee
approximates the true least-squares
Fourier coefficients well apart from
an $O(1)$ number of them lying at extreme frequencies near $m=\pm N/2$.
For $R>\rho$ these Fourier coefficients grow exponentially,
so dropping an $O(1)$
factor we may consider only the largest coefficient's contribution
to the $l^2$-norm, namely that at $m=N/2$.
For $R<\rho$ there is no growth hence the norm is
dominated by the coefficients of size $O(1)$ at $m\approx 0$.
Therefore
we have the order-of-magnitude estimate at a given $N$,
\be
|\mbf{\alpha}| \approx \mbox{max}\left[\;\frac{1}{\sqrt{N}} \,
\frac{|\hat{v}(N/2)|}{|\hat{s}(N/2)|}\; ,\; 1\; \right],
\label{eq:nap}
\ee
where the maximum-value operation combines the two cases.
Similarly, since the boundary data coefficients die exponentially,
following Remark~\ref{rmk:interp} and ignoring an $O(1)$ factor
we may suppose
\be
t \approx |\hat{v}(N/2)|.
\label{eq:tap}
\ee
We define $N_0$ to be the $N$ at which convergence stops,
for this we use the round-off error consideration
$t/|\mbf{\alpha}| \approx \emach$ discussed above.
In the case $R<\rho$ this implies that convergence halts
when $t$ reaches of order $\emach$; as observed in
Fig.~\ref{fig:lowk}a,b.
However for $R>\rho$, combining \eqref{eq:nap} and \eqref{eq:tap} we get an
implicit equation for $N_0$,
\be
\sqrt{N_0}|\hat{s}(N_0/2)| \;\approx\; \emach
\qquad \mbox{(criterion for halting of convergence, for $R>\rho$)}~.
\label{eq:halt}
\ee
The minimum achievable boundary error is then
given by \eqref{eq:tap} with the substitution $N=N_0$.
As an illustration,
using the (Laplace) asymptotic form \eqref{eq:lapas} for the eigenvalues
approximately predicts (dropping algebraic factors) that
$N_0 \approx 2 \ln (1/\emach) / \ln R$.
For the parameters of Fig.~\ref{fig:lowk}c this gives
$N_0 \approx 180$, then using \eqref{eq:decay} with $C=1$ gives
$t_0 \approx 10^{-4}$ which, given the heuristic nature of our model,
agree well with the observed behavior.

We now briefly discuss the case of high wavenumber.
In Fig.~\ref{fig:sm} we saw that the (Laplace) asymptotic form
\eqref{eq:lapas}
is not useful for predicting relevant eigenvalues at high $k$.
We may derive (Appendix \ref{app:bes}) the asymptotic
\be
\hat{s}(m) \sim \frac{1}{2|m|} R^{-|m|} e^{k^2(R^2-1)/4m},
\qquad |m|\to\infty,
\label{eq:smser}
\ee
which Fig.~\ref{fig:sm}a),~b) shows is a much improved approximation,
but still not useful for the relevant
eigenvalues at $k=500$ (or beyond).
Therefore we use the WKBJ method to
derive (see Appendix \ref{app:bes})
a uniform approximation for the eigenvalue magnitudes,
defining $a^2=m^2-\frac{1}{4}$,
\be
|\hat{s}(m)| \approx \left\{\begin{array}{ll}
\left[(k^2 - a^2)(k^2R^2-a^2)\right]^{-1/4}, & m<k\\
\frac{1}{2}\left[(a^2-k^2)(k^2R^2-a^2)\right]^{-1/4} e^{I_a(k)}, & k<m<kR\\
\frac{1}{2}\left[(a^2-k^2)(a^2-k^2R^2)\right]^{-1/4} e^{I_a(k)-I_a(kR)},& m>kR
\end{array}\right.
\label{eq:unif}
\ee
where
\be
I_a(x):=\sqrt{a^2-x^2} - a \ln[(a+\sqrt{a^2-x^2})/x].
\label{eq:Iax}
\ee
More precisely this is an estimate of the {\em amplitude} in
oscillatory region ($m<k$) of $J_m$, and the absolute value in
the evanescent region (note $H_m^{(1)}$ can also be complex oscillatory
but its magnitude never is).
In the oscillatory region
individual $\hat{s}(m)$ values cannot be predicted: rather they
are distributed in the range $[-1,1]$ times the approximate
amplitude \eqref{eq:unif}.
Fig.~\ref{fig:sm} shows this is a highly accurate asymptotic form
in all regions apart from the two turning-points ($J_m(k)$ is at its turning
point for $m\approx k$ whereas for $H_m^{(1)}(kR)$
this occurs at $m\approx kR$).
The estimate has algebraic singularities at these two turning-points, but
they are weak enough that it is still useful.

Finally, we compare numerical convergence results against this model at
high wavenumber.
Fig.~\ref{fig:highk} shows convergence and coefficient norm
at $k=500$ (about 170 wavelengths across the domain)
for boundary data \eqref{eq:vfund} deriving from an exterior fundamental
solution.
Its boundary data Fourier coefficients $\hat{v}(m)$
are given by the same formula \eqref{eq:sm} as the MFS eigenvalues
(and hence the same approximation \eqref{eq:unif})
but with the substitution $\rho$ for $R$.
To compute the curves shown as `predicted',
we used this approximation in \eqref{eq:nap},
and \eqref{eq:tap} to predict the error norm.
It is clear that, up to the point when convergence halts,
the predictions for both error norm and coefficient norm
are very close to observations (the largest deviations being
spikes due to algebraic singularities discussed above;
in a) these are at $N=1000$ and 1100).

A crucial common feature is that no convergence happens until $N=2k$,
since $\hat{s}(m)$ remains large for $|m|<k$ (see Fig.~\ref{fig:sm}b).
One interpretation of this is that $2k$, corresponding to
2 degrees of freedom per wavelength on the perimeter, is the
Nyquist sampling frequency for $k$-bandlimited functions on $\pO$;
in physics this is known as the the `semiclassical basis size' \cite{que}.
In panel a) of the figure, $\rho=1.1$ so the singularity is
$k(\rho-1)/2\pi \approx 8$ wavelengths from the boundary.
In this case convergence is rapid, dropping ten orders of magnitude
between $N=1000$ and $N=1150$. Convergence then halts (compare \eqref{eq:halt}
which predicts $N_0\approx 1175$ and $t_0\approx3\times 10^{-14}$).
The number of boundary quadrature points was
$M=1500$ in a) (only 3 points per wavelength).
This can be chosen to be so small
since boundary functions $v$ and $u^{(N)}|_{\pO}$
have exponentially-decaying Fourier coefficients beyond frequency $2k$,
giving spectral convergence.
In panel b) $\rho=1.01$ and $R=1.05$, giving both slower convergence
and growth in coefficient norm.
The predictions $N_0 \approx 1708$ and $t_0 \approx 8\times 10^{-7}$
are again reasonably close to observations.

How can $R$ best be chosen to achieve the lowest boundary error
for a given high wavenumber $k$, and $\rho$? We use the above model to compute
$N_0$ and hence $t_0$ for a variety of $\rho$ and $R$ at $k=500$,
in Fig.~\ref{fig:Rsweep}.
Here the smallest $\rho=1.002$ corresponds to a singularity 0.16 wavelengths
from the boundary.
The conclusion is that $t_0 \approx \emach$ appears to be always
achievable as $R$ tends to $\rho$ from above, as expected from
Remark~\ref{rmk:bestR};
however, the basis size required to do this diverges as $\rho\to1^+$.
b) also shows that there is a limiting basis size of about $N\approx1180$
(not much larger than $2k$)
for which arbitrarily large $R$ may be used,
but with this choice $t_0$ becomes $O(1)$, hence not useful, as $\rho\to1^+$.
In conclusion, we may state that
in the high-wavenumber limit, if the nearest singularity in the
boundary data is at least a few wavelengths away, then
both the basis size $N$ and the number of quadrature points $M$
can approach 2 per wavelength while achieving an error
close to machine precision.




\section{The MFS on analytic domains}
\label{sec:mfsanalytic}

In this section we present results for the MFS on arbitrary analytic
domains. On the circle we have shown that the MFS coefficients
start growing exponentially if the radius $R$ of the charge points
becomes larger than the distance $\rho$ of the singularity. In this
section we demonstrate that also on general analytic domains the
position of the charge points relative to the singularities of the
analytic continuation is crucial for the accuracy and numerical
stability of the MFS.

\subsection{Analytic continuation of solutions}
\label{sec:continuation}

The question of analytic continuation is to
find a domain $\cOmega\supset\Omega$ and a function $\cu$ such that
$\Delta \cu+\lambda \cu=0$ in $\cOmega$ and $\cu|_{\Omega}=u$. Since
solutions of the Helmholtz equation are real analytic it follows
immediately that $\cu$ is unique.

A classical result of analytic continuation is reflection on a
straight arc $\Gamma$, on which $u$ satisfies $u|_\Gamma=0$. Without
restriction let $\Gamma$ be a subset of $\{iy:~y\in\mathbb{R}\}$.
Then $u$ can be continued across $\Gamma$ by setting
$u(-x,y):=-u(x,y)$ (see also \cite{CoHi53}). In \cite{Ga53}
Garabedian extended these results to the case that $\Gamma$ is an
arbitrary analytic arc for which $u_\Gamma=0$. More general
reflection principles for linear elliptic PDEs of the type $\Delta
u+a(x,y)u_x+b(x,y)u_y+c(x,y)=0$, where $a(x,y)$, $b(x,y)$ and
$c(x,y)$ are real analytic functions were treated by Lewy in
\cite{Le59}. He stated his results for arbitrary Dirichlet, Neumann
and mixed boundary conditions but restricted $\Gamma$ to be a
straight line. Representations of the analytic continuation for the
case that $\Gamma$ is not a straight line were given by Millar in
\cite{Mi80}. In \cite{Mi86} he discussed more in detail the analytic
continuation of solutions of the Helmholtz equation.

Millar shows that there are two possible sources for singularities of
the analytic continuation $\cu$ of $u$. The first one comes from
singularities of the analatic continuation of the boundary data
$f$. The second possible source of singularities is introduced by the
shape of $\partial\Omega$. Let $Z(s)=x(s)+iy(s)$ be a parameterization
of $\partial\Omega$, where $s\in[0,2\pi]$. Assume that $x(s)$ and
$y(s)$ are real analytic and that $|Z'(s)|\neq 0$ in $[0,2\pi]$. Then
there exists a complex neighborhood of $[0,2\pi]$, in which $Z(s)$ is
holomorphic and invertible. We denote its inverse by $S(z)$ and define
the Schwarz function
$$
G(z):=\bar{Z}(S(z))=\overline{Z(\overline{S(z)})}.
$$
Millar showed that except for special cases the
singularities of $G(z)$ outside $\Omega$ are also singularities of the analytic
continuation $\cu$ of $u$.

\begin{figure}
\center
\includegraphics[width=12cm]{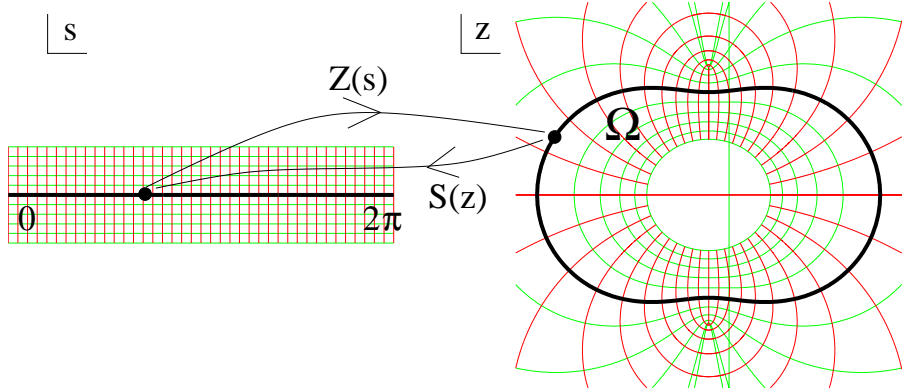}
\caption{Illustration of map $Z(s)$ and its inverse $S(z)$ defining
a boundary curve $\pO$.
The Schwarz function involves composition of $Z(s)$ and $S(z)$.}
\label{fig:map}
\end{figure}

The Schwarz function has been studied in \cite{Da74}. It is
independent of the parameterization of $\partial\Omega$ and has an
interpretation in terms of reflection principles on analytic
arcs. Assume that $z_1$ is a point close to $\partial\Omega$. Then its
reflection on $\partial\Omega$ can be obtained by the following steps
(see Figure \ref{fig:map}).
\begin{enumerate}
 \item Compute $t_1=S(z_1)$.
\item Reflect $t_1$ on the real line to obtain the point $t_2:=\overline{t_1}$.
\item The reflection $z_2$ of $z_1$ at $\partial\Omega$ is now obtained as
\begin{equation}
\label{eq:zreflection}
z_2=Z(t_2)=Z(\overline{t_1})=Z(\overline{S(z_1)})=\overline{G(z_1)}.
\end{equation}
\end{enumerate}

\begin{figure} 
\center
 \includegraphics[width=13cm]{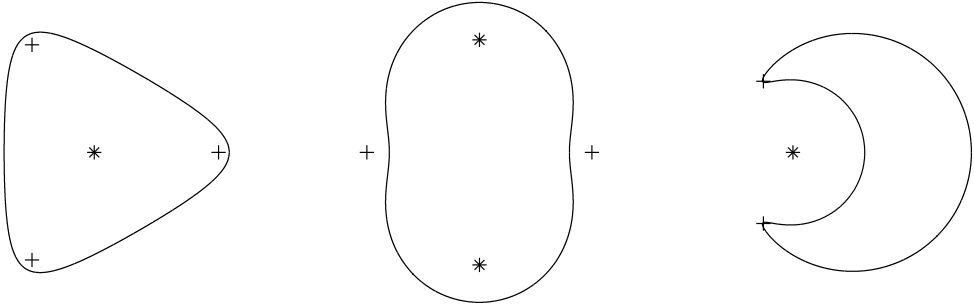}
\caption{
Domains a) rounded triangle, b) inverted ellipse, c) crescent.
Branch singularities of the Schwarz function are denoted by '+' and
pole type singularites by '*'.
Note for the crescent, the '+' signs are inside $\Omega$
but very close to the boundary.
}
\label{fig:domains}
\end{figure}

Fig. \ref{fig:domains} shows the singularities of the Schwarz function on
three different domains, a rounded triangle,
an inverted ellipse and a crescent. The
domains are defined with default values for the
parameters $a_1$ through $a_4$, as follows.

\begin{tabular}{lc}
Rounded triangle: &$Z_T(s)=e^{is}+a_1e^{-2is},$ \qquad $a_1 = 0.3$\\
Inverted ellipse: & $Z_{IE}(s)=\frac{e^{is}}{1+a_2e^{2is}}$,
\qquad $a_2 = 0.25$\\
Crescent: & $Z_C(s)=e^{is}-\frac{a_3}{e^{is}+a_4}$, \qquad $a_3=0.1, a_4=0.9$
\label{t:shapes}
\end{tabular}

Branch type singularities are denoted by '+' and pole type
singularities by '*' in Figure \ref{fig:domains}. The
branch singularities in all three domains are of square root type (see
\cite{Mi86} for an analysis of the branch behavior of $G(z)$).
The crescent has exterior singularity of pole type at $z=-1/\overline{a_4}$.
For the
interior Helmholtz problem only the exterior singularities of $G$ are
important since these are points where, for generic boundary data, the
analytic continuation $\cu$ of $u$ becomes singular. Conversely if we
had an exterior Helmholtz problem then the interior singularities
would determine the singularities of the analytic continuation.

\subsection{Using exterior conformal map to place the charge points}
\label{sec:conf}

A natural generalization of the MFS on the unit disk to general
analytic domains can be defined in terms of the conformal map from the
exterior of the unit disk to the exterior of the domain.
This was investigated in the Laplace BVP case by
Katsurada \cite{Ka94}.

Let $\Omega$ be a simply connected domain with analytic boundary
$\partial\Omega$. We can parameterize $\partial\Omega$ using the
exterior conformal mapping function
$$
z =
\Psi(w)=cw+c_0+\frac{c_{-1}}{w}+\frac{c_{-2}}{w^2}+\dots,\quad c>0
$$
which maps the exterior $D_1:=\{w:~|w|>1\}$ of the unit disk to the exterior of $\Omega$.
The quantity $c$ is called the capacity of
$\Omega$. We denote the inverse map by $w = \Phi(z)$.  Since
$\partial\Omega$ is analytic $\Psi(w)$ can be analytically continued
to a domain $D_r:=\{w:~|w|>r\}$ for some $0<r<1$.
We denote the conformal radius of a point
$z\in\mathbb{C}\backslash\Omega$ by $\rho_z:=|\Phi(z)|$.

In the notation of the previous section
we may write this parametrization as $Z(s) = \Psi(e^{is})$
since the unit disc is parametrized by $w=e^{is}$.
Using that the reflection of a point $z$ on the unit circle is given by
$z'=\frac{1}{\overline{z}}$ the Schwarz function may now be written
$G(z)=\overline{\Psi(1 / \overline{\Phi(z)})}$; it follows that it is
analytic in $\{z\in\mathbb{C}:1<\rho_z<\frac{1}{r}\}$.

In the unit disk case we placed the MFS points equally distributed on
a curve with radius $R$.
For general analytic domains we now place the
points on a curve $\Gamma_R := \{z: \rho_z = R\}$
with constant conformal radius $\rho_z=R$. On this curve we distribute
the points equally spaced in conformal angle, that is
\be
\by_j:=\Psi(e^{2\pi ij/N}), \qquad j=1,\dots,N
\label{eq:yextconf}
\ee
If $\Omega$ is the unit disk this definition coincides with
that of Sec.~\ref{sec:unitdisc}. Replacing the disk radii $R$ and
$\rho$ in Theorem \ref{thm:t} by the corresponding conformal radii we
obtain the following conjecture for the rate of convergence of the MFS
for Helmholtz problems on general analytic domains.

\begin{cnj} 
Let $t$  be the error of the MFS as defined in \eqref{eq:t} by placing
the MFS points equally distributed in conformal angle at a conformal
distance $R$ around $\Omega$. Let $\rho>1$ be the conformal radius
of the closest (in the sense of conformal radius) singularity of
the analytic continuation of $u$.  Then
\be
t \; \le \; \left\{\begin{array}{ll}
C\rho^{-N/2},& \rho < R^2,\\
CR^{-N},& \rho > R^2,
\end{array}\right.
\ee
\label{cnj:t} 
where $C$ is a constant that may depend on $\Omega$, $k$, $R$ and $v$, but not
$N$. Furthermore, if $u$ continues to an entire function,
the latter case holds for any $R>1$.
\end{cnj}
\begin{rmk} 
The first case of this conjecture was proved in the Laplace case by
Katsurada in \cite{Ka94} under additional restrictions on the analytic
continuation of $\Psi$ into the unit disk.
In numerical studies we have observed that these
conditions are not necessary to achieve the given convergence rates,
so do not include them in our conjecture (compare also Remark 3.2 of \cite{Ka94}).
We do not state a conjecture for the case $\rho=R^2$ since
numerically it cannot be established if for general domains
the same algebraic factor is needed as for the disk.
\end{rmk} 

\begin{figure} 
\center
\mbox{
\includegraphics[width=5in]{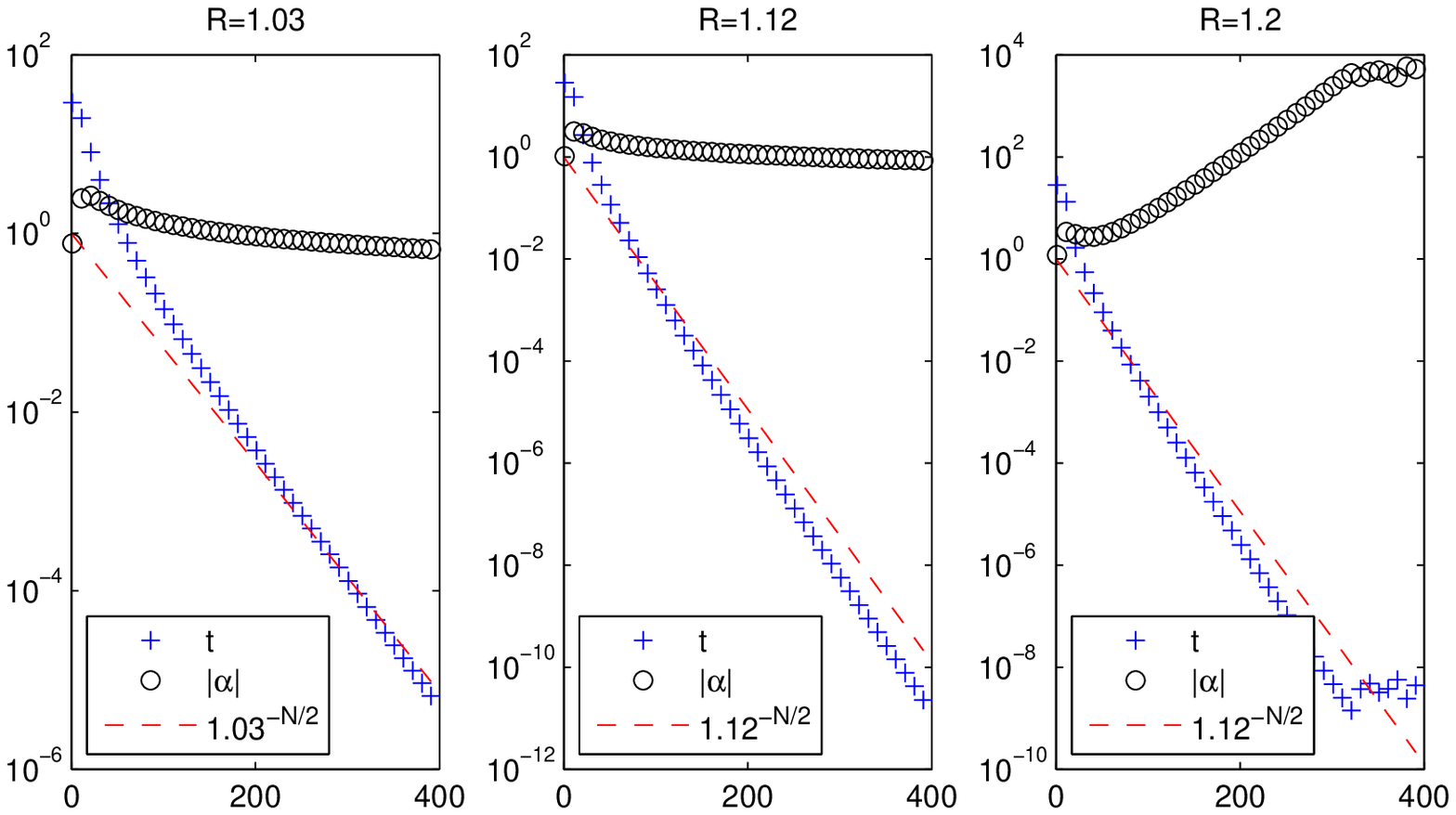}
\includegraphics[width=1.5in]{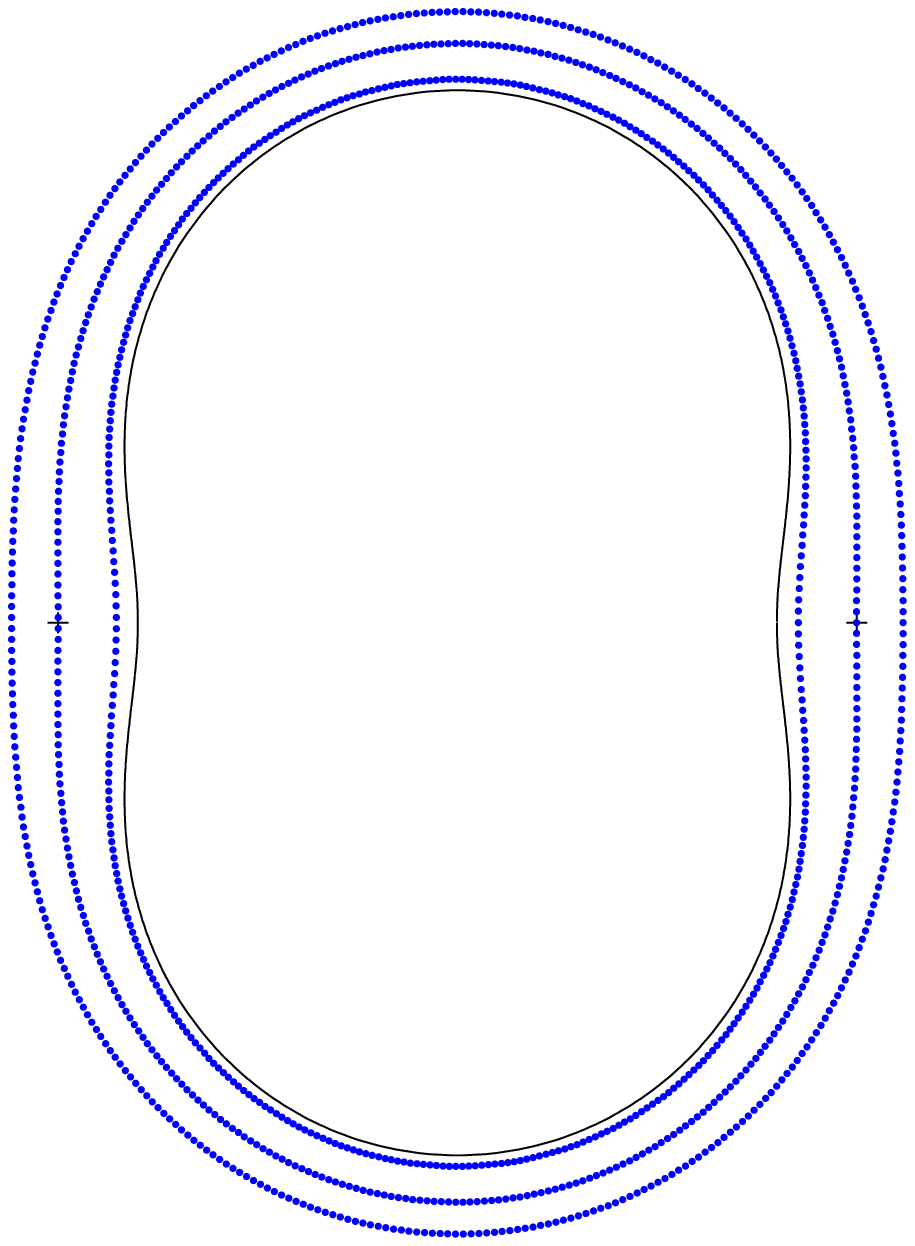}
}
\caption{Left three plots: Estimated (dashed lines) and observed ('+') rates
  of convergence of the MFS on the inverted ellipse for charge
  points with conformal radii $R=1.03$, $R=1.12$ and $R=1.2$. The
  corresponding coefficient norms $|\bmal|$ are denoted by 'o'.
Right
  plot: Domain and positions of the charge points (shown as dots)
for the above three
conformal radii, with branch-type singularity locations ('+' signs).}
\label{fig:conformalconv}
\end{figure}

In Figure \ref{fig:conformalconv} we plot the observed error $t$ for
the MFS on the inverted ellipse of Figure \ref{fig:domains}, for
wavenumber $k=5$ and constant boundary condition $v\equiv1$. The estimated
rates from Conjecture \ref{cnj:t} are denoted by dashed lines.
The three plots correspond to MFS
points placed at the conformal distances $R=1.03$, $R=1.12$ and
$R=1.2$. The conformal radius $R=1.12$ is also the approximate
conformal radius $\rho$ of the two singularities. The corresponding
MFS curves are shown on the right of Fig.~\ref{fig:conformalconv}.
The estimated convergence rates are in all
three cases in good agreement with the observed error $t$.

In the disk case we can observe exponential coefficient growth once
the radius $R$ of the charge points is larger than the radius $\rho$
of the singularity of the analytic continuation of $u$ (see
Thm. \ref{thm:a}). For general analytic domains we observe a
similar behavior. To demonstrate this we plotted in Figure
\ref{fig:conformalconv} also the norm $|\bmal|$ of the MFS
coefficients. As long as the conformal radius $R$ of the MFS points is
smaller or equal to the conformal radius $\rho$ of the singularities
we do not observe any growth of $|\bmal|$ for growing $N$ (first two
plots of Fig. \ref{fig:conformalconv}). In the third plot we have
$R>\rho$ and $|\bmal|$ grows exponentially for growing $N$.
It is instructive to compare this figure panel by panel against
Fig.~\ref{fig:lowk}.

This leads to the following conjecture, which mirrors Theorem
\ref{thm:n} and \ref{thm:a}.
\begin{cnj} 
Let the MFS charge points be chosen equally spaced in conformal angle on
$\Gamma_R:=\{z:~\rho_z=R\}$, the curve of all points with
given conformal radius $R>1$. Let $\rho$ be the conformal distance
of the closest singularity of the analytic continuation of $u$ and let
$\bmal$ be the vector of coefficients of the MFS basis functions that
minimizes $t$.
We have
$$
|\bmal| \ge C\gamma^{N}
$$
for some $\gamma>1$, and a constant $C$ that does not depend on $N$, if
and only if $R>\rho$.
\label{cnj:a} 
\end{cnj}
In numerical experiments we observed also for other types
of MFS curves that there is only coefficient growth if the curve
encloses a singularity of the analytic continuation. Hence, a more
general conjecture can be stated (similar to results known for the
scattering case~\cite{Ky96}).
\begin{cnj} 
  Let $\Gamma$ be any Jordan curve enclosing $\overline{\Omega}$,
  with $\text{dist}(\Gamma,\pO)>0$,
  on which MFS charge points are chosen asymptotically densely.
  Then the
  coefficient norm $|\bmal|$ that minimizes $t$ 
grows asymptotically exponentially as
  $N\rightarrow\infty$ if and only if $\Gamma$ encloses a
  singularity of the analytic continuation of $u$.
\label{cnj:agen}
\end{cnj} 

\begin{figure} 
\center
\includegraphics[width=0.9\textwidth]{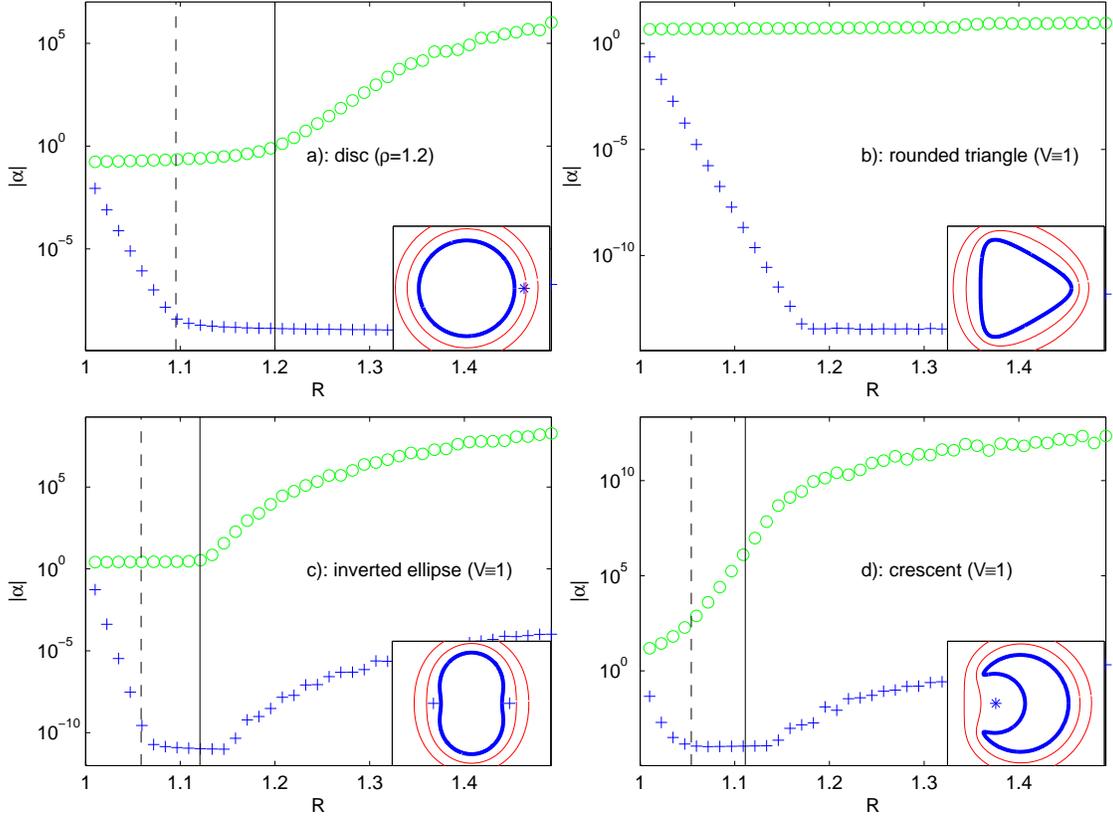}
\caption{
Growth of coefficient norm $|\bmal|$ (circles),
and least-squares approximation error $t$ ('+' signs), for a) the unit disc
with boundary data \eqref{eq:vfund}, and
b-d) the three other shapes of Fig.~\ref{fig:domains}
with constant boundary data $v\equiv 1$.
MFS points are chosen on curves with constant conformal radius $R$.
The insets show domains, exterior singularities
using same notation as Fig.~\ref{fig:domains},
and curves with conformal radius $R=1.24$ and $1.49$.
The parameters are $k=5$. $N=200$ for a-b, $N=400$ for c-d. 
The number of quadrature points $M$ on $\partial\Omega$ is
chosen sufficiently large throughout.
}
\label{fig:coeffgrow}
\end{figure}

In Figure \ref{fig:coeffgrow} we show the
coefficient norm $|\bmal|$
and the approximation error $t$ for a growing conformal
distance $R$ of the MFS source points and fixed number $N$ in four different cases: the
unit disk and the three domains from Figure~\ref{fig:domains}. In all
cases we have used $k=5$. For the disk the boundary data is given by
\eqref{eq:vfund} with singularity location $\rho=1.2$,
and in the other three cases by $v(z)\equiv1$
(recall that here the Schwarz function introduces singularities in $u$).
The vertical solid lines denote the conformal
radius $\rho$ of the singularities of the analytic continuation of $u$ and
the vertical dashed lines denote the square root $\rho^{1/2}$.
Since the Schwarz function for the rounded triangle does not have any
singularities in the exterior of
the domain the solution $u$ can be analytically continued to an entire
function.

For the disk, the inverted ellipse and the crescent the error $t$ does not
decrease further once $R$ passes the dashed line. This can be expected
from Conjecture \ref{cnj:t} since the upper bound on the error $t$
does not decrease any more for fixed $N$ and $R>\rho^{1/2}$.

For these three domains we can also observe exponential coefficient growth of
$|\bmal|$ for fixed $N$ when $R>\rho$. For the crescent this exponential growth
already starts earlier. However, this is not a contradition to
Conjecture \ref{cnj:a}. The conjecture treats the case of fixed $R$
and $N\rightarrow\infty$. This does not exclude the existence of
transient growth effects for $R<\rho$. An explanation for these
transient effects in the crescent case is that close to the pole-like
singularity of the
Schwarz function we need a very high number $N$ of basis functions to
sufficiently resolve a highly-oscillatory Helmholtz field.

Another interesting special case is the rounded triangle. Since the
analytic continuation of $u$ is an entire function, by Conjecture
\ref{cnj:a} we do not expect any exponential growth of
$|\bmal|$. Indeed, Fig.~\ref{fig:coeffgrow}b shows that
$|\bmal|$ stays virtually constant as $R$ increases.

In this section we have demonstrated with several numerical
experiments that the behavior of the MFS for Helmholtz problems on
general analytic domains is similar to the unit disk case if we choose
the MFS points on curves with constant conformal radius $R$. We stated
two conjectures about the approximation error $t$ and the coefficient
norm $|\bmal|$ as $N\rightarrow\infty$. From our numerical
experiments and under the conditions that the conjectures hold it
follows that the optimal conformal radius $R$ for the position of the
singularities is given by $R=\rho^{1/2}$, where $\rho$ is the
conformal radius of the closest singularity. This ensures a maximum
rate of convergence for $t$ while keeping $|\bmal|$ from growing
exponentially
(as in Remark~\ref{rmk:bestR}).
Furthermore, in view of the results in the crescent
case of Figure \ref{fig:coeffgrow} it seems advisable to choose $R$
not too close to $\rho$ if the singularity is determined by a pole in
the Schwarz function in order to avoid large transient coefficient
growth.

\begin{figure} 
\center
\includegraphics[width=0.7\textwidth]{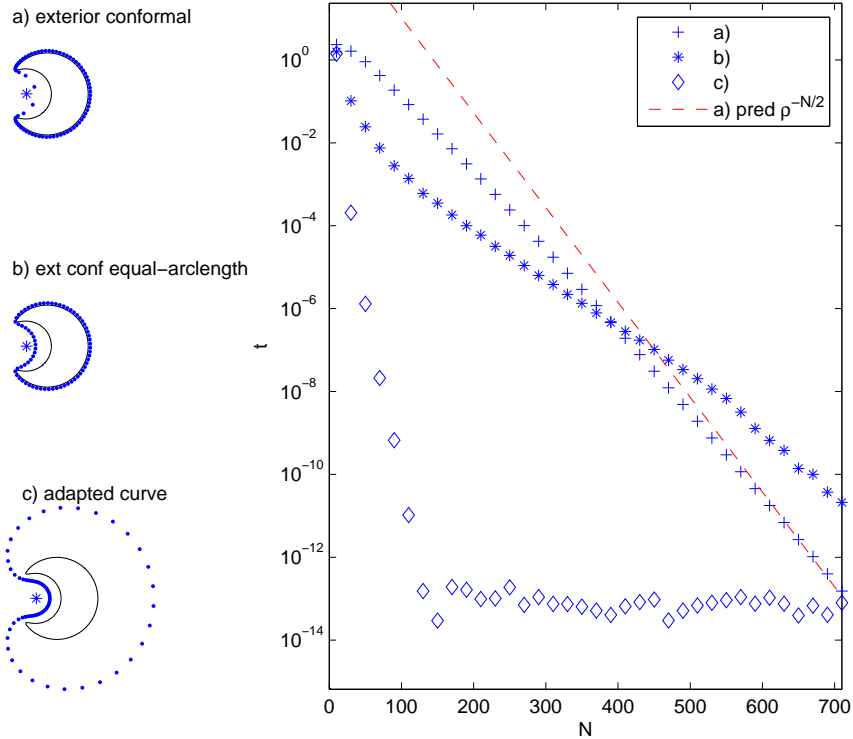}
\caption{\label{f:beatconf}
Convergence rate for the crescent domain (with $a_3=0.1, a_4=0.9$)
illustrating the concavity problem.
The wavenumber is $k=3$ and constant boundary data $v\equiv1$.
For a) and b) MFS charge points are placed on an exterior conformal curve
at $R = \sqrt{\rho}$ where $\rho = 1/a_4$ is the singularity conformal
distance. The point spacing is equal in a) conformal angle, and b) arclength.
c) Adapted curve and spacing given by \eqref{eq:adapt}, see
Sec~\ref{sec:adapt}.
The dashed line shows the predicted convergence rate for a), $a_4^N$.
}
\end{figure} 

\begin{figure} 
\center
\includegraphics[width=\textwidth]{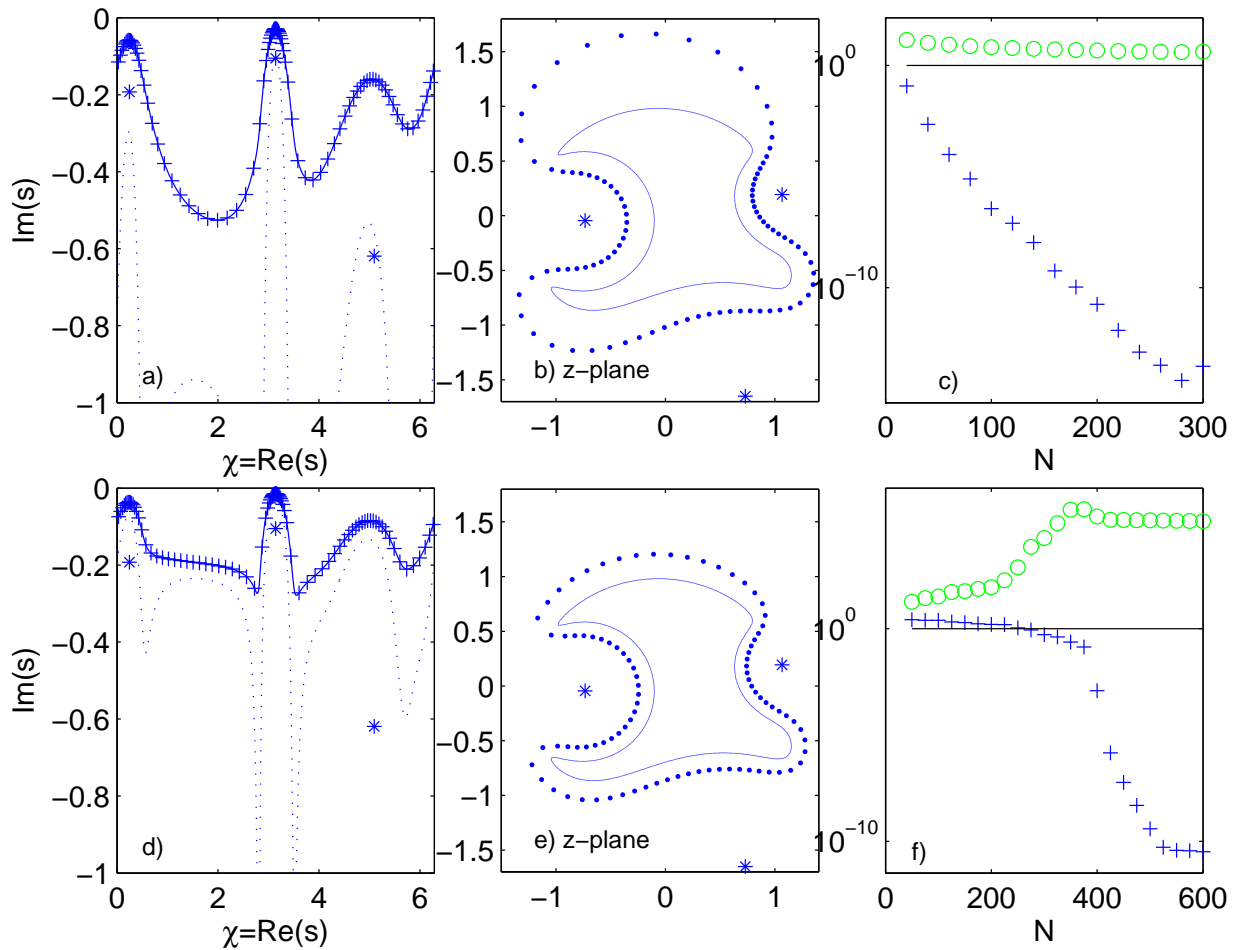}
\caption{\label{f:splane}
Generalized crescent analytic domain parametrized by \eqref{eq:gencres},
showing performance of adaptive charge points of Sec.~\ref{sec:adapt}.
The top row is at $k=3$, the bottom row at $k=100$.
a) and d) show the curve $y(\chi)$ and the locations $(\chi_j,y(\chi_j))$
in the $s$-plane, the singularities $s_\sigma$ ('*' symbols),
and the distance-limiting function $|Z'(\chi)|/\dmax$ (dotted line).
b) and e) show the charge locations (for clarity, $N=90$ has been
used in both cases).
c) and f) shows error norm convergence $t$ ('+' symbols) and
coefficient norm $|\bmal|$ ('o' symbols).
}
\end{figure} 

\begin{figure} 
\center
\includegraphics[width=0.6\textwidth]{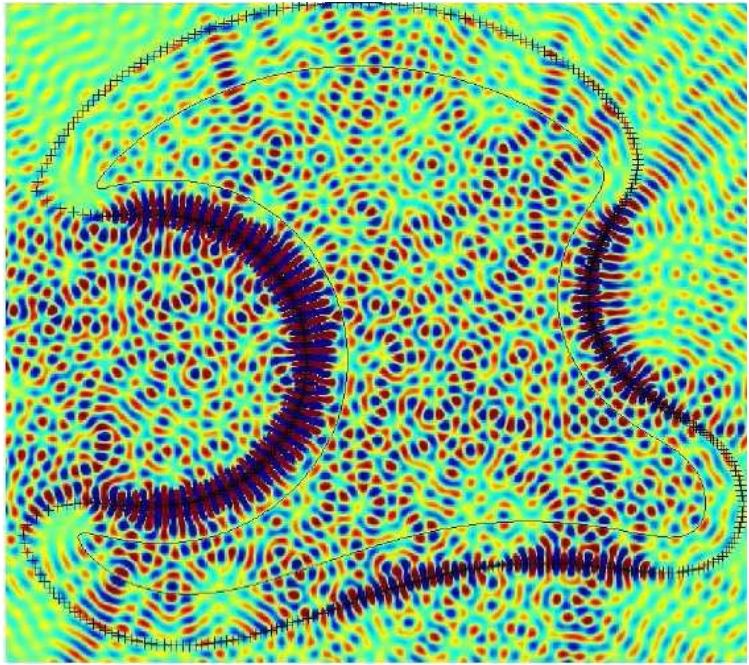}
\caption{\label{f:gencresug}
Interior Helmholtz BVP solution $u$ for
the same domain as Fig.~\ref{f:splane}, at $k=100$,
with boundary data $v\equiv 1$,
using as basis size of $N=525$, and $M=1000$ boundary
points.
Error norm is $t=5\times 10^{-11}$.
The boundary (thin solid line) and charge points ('+' symbols) are shown,
and the field $u$ from \eqref{eq:mfs} is shown both outside and inside
$\Omega$.
CPU time was 2.9 s to compute the coefficient vector $\bmal$,
and 5.4 mins to evaluate the solution $u$
shown ($2.4\times 10^5$ points on a grid of spacing 0.005).
}
\end{figure} 

\begin{figure} 
\center
\includegraphics[width=\textwidth]{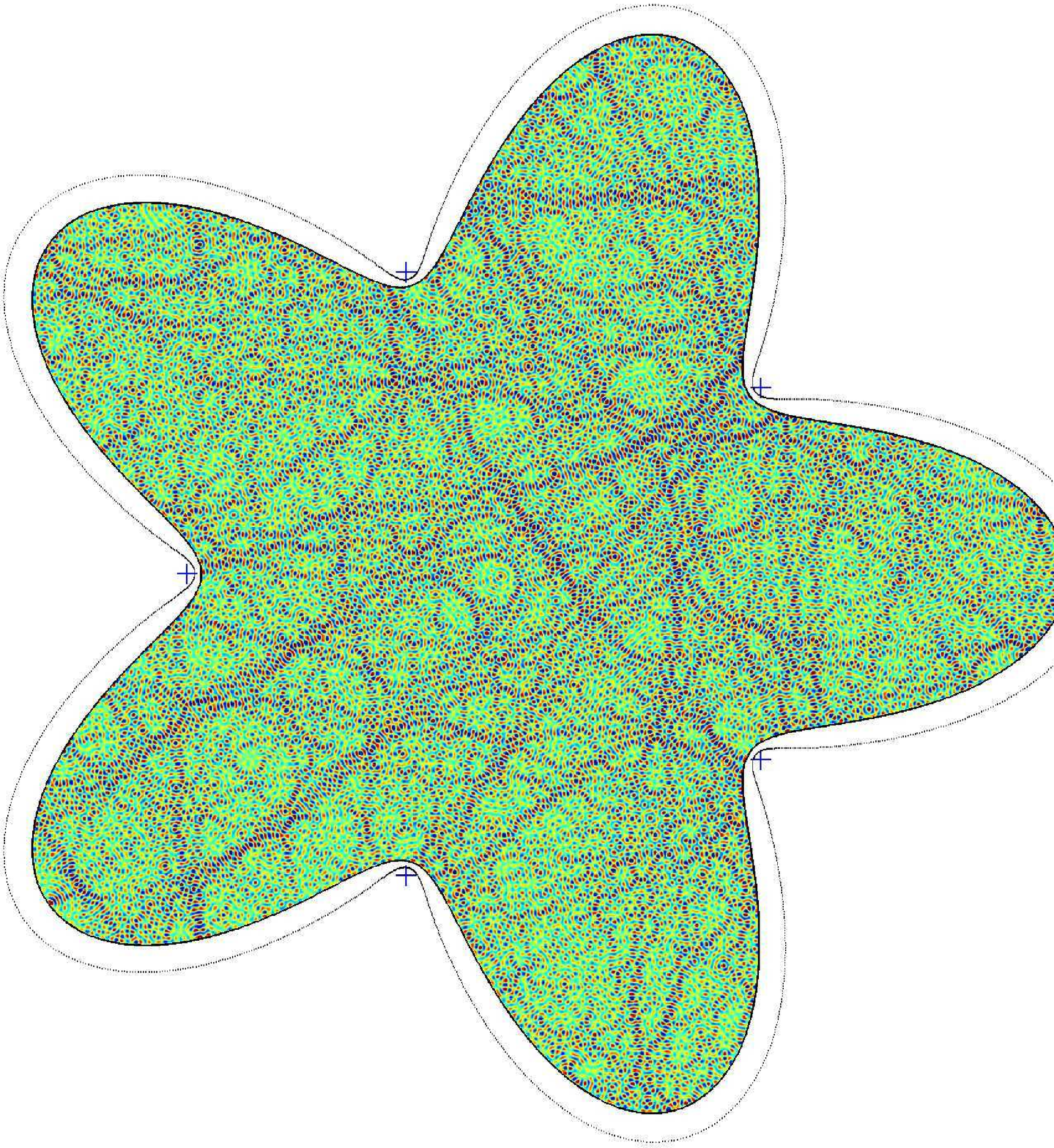}
\caption{\label{f:pentafoilug}
Interior Helmholtz BVP solution $u$ for
the analytic 5-foil domain given
by the radial function $r(\theta) = 1+ 0.3 \cos(5\theta)$
at wavenumber $k=400$,
with boundary data given by $v(z) = \mbox{Re}\,(z-\rho)^{-1}$
with $\rho=1+0.5i$.
Error norm is $t=2\times 10^{-11}$.
Basis size was $N=2000$, and $M=2000$ boundary
points. Charge points are shown by dots and Schwarz branch singularities
by '+' symbols.
CPU time was 26 s to compute the coefficient vector $\bmal$,
and 1.05 hr to evaluate the interior solution $u$
shown ($8.2\times 10^5$ points on a grid of spacing 0.002).
}
\end{figure} 

\subsection{Using a singularity-adapted curve to place the charge points}
\label{sec:adapt}

The above exterior conformal method has the following problem:
in any concave parts of $\Omega$ the exterior conformal map $\Psi$
has a very (in fact, exponentially) large gradient.
This well-known property of conformal maps is related
to the so-called crowding problem.
This has two consequences for concave regions:
the spacing of charge points according to \eqref{eq:yextconf}
becomes very large, and
as $R$ is increased from 1, the curve $\Gamma_R$ moves away from
$\pO$ very rapidly.
Both these effects are illustrated by the
MFS charge curve for the crescent in Fig.~\ref{f:beatconf}a.

This means that Schwarz function
singularities which are a moderate distance from a concave part of
$\pO$ may actually have a conformal radius extremely close to 1.
The net result is that if the coefficient growth of Conj.~\ref{cnj:a}
is to be prevented, $R$ must be very close to 1, hence
by Conj.~\ref{cnj:t} the convergence rate is necessarily very poor.
For example, in Fig.~\ref{f:beatconf},
the conformal radius of the pole in the
Schwarz function is only $\rho = 1/a_4 \approx 1.11$, and
the observed rate for case a) approaches the predicted $\rho^{-N/2}$
(dashed line in the figure).
One can attempt to fix the problem of the large point spacing
by retaining the same MFS curve $\Gamma_R$ but choosing charge points
equally spaced in arc-length, as illustrated in  Fig.~\ref{f:beatconf}b.
Despite an initial improvement for small $N$, the
asymptotic convergence rate turns out to be no better than in case a),
and is believed to be the same
(for errors $t<10^{-6}$ it performs worse, we believe due to
a lack of point density near the `spiked' parts of the crescent).

In Fig.~\ref{f:beatconf}c we show a different `adaptive' choice of MFS curve
and point spacing which
clusters charge points near the singularity but spreads them out
(while taking them further away from $\pO$) away from the singularity.
This gives a convergence rate over 5 times faster than the exterior
conformal curve, for instance
an error of $t\approx 10^{-13}$ is reached for $N=140$ as opposed to $N=730$
for case a).
This simple adaptive method is based on the idea of replacing
the exterior conformal map by an annular one, as follows.

From the discussion in Sec.~\ref{sec:continuation}, the map $Z(s)$
defines an annular conformal map which is
one-to-one for $s$ in some strip around $[0,2\pi]$, as in Fig.~\ref{fig:map}.
The external singularities control convergence rate;
we label them by $\sigma = 1,\ldots, P$, where $P$
is the number of singularities.
They have $s$-plane locations $s_\sigma = \chi_\sigma-i\tau_\sigma$
with $\tau_\sigma > 0$,
Their minimum distance to the real axis is $\tau:=\min_\sigma \tau_\sigma$.
Katsurada \etal~\cite{Ka96} have discussed using such an annular map to place
charge points for the Laplace BVP, according to
$\by_j:=Z(2\pi j/N-i \log R)$ for some $R>1$;
note this is the annular map equivalent of \eqref{eq:yextconf}.
A related annular map method has been tried
in the scattering case \cite{Kark01}.
According to Conj.~\ref{cnj:agen}
in order to prevent coefficient growth in this case
one would need to choose $\log R < \tau$.
This may be a severe restriction: for example one may
check that the crescent formula $Z_C(s)$ in Table~\ref{t:shapes}
is identical to this domain's exterior conformal map $\Psi$,
thus the concavity effect causes $\tau$ to be very small.

However, there are an infinite family of annular maps $Z(s)$ which
parametrize the {\em same} boundary $\pO$, and these may have
differing Schwarz singularity locations in the $s$-plane.
Given an analytic
domain $\Omega$ one is free to choose between such parametrizations.
Ideally, our goal is to choose one with as large a $\tau$ as possible,
to achieve a high convergence rate.
However, we find it convenient to retain the given parametrization $Z(s)$,
and instead build a charge curve in the $s$-plane which
no longer has constant imaginary part, and which captures
the spirit of such a reparametrization.
Our curve is given by $s = \chi - iy(\chi)$, where $\chi \in [0,2\pi]$
is the real part of $s$, and the function $y$ is given by
\be
[y(\chi)]^{-1} = \left(\frac{\dmax}{|Z'(\chi)|}\right)^{-1} +
\sum_{\sigma=1}^{P}\left[\gamma \tau_\sigma +
  \beta\frac{1-\cos(\chi-\chi_\sigma)}{\tau_\sigma}\right]^{-1}
\label{eq:adapt}
\ee
where parameter values performing well in most domains
are $\beta = 0.7$ (interpreted as a curvature factor),
$\gamma = 0.4$ (expressing the curve's fractional
distance to each singularity), and
\be
\dmax = \mbox{max}[1, \;25/k]
\label{eq:dmax}
\ee
is interpreted as the
maximum allowable distance of the curve from $\pO$.
Roughly speaking, \eqref{eq:adapt} has the effect of bringing the
curve close to $\pO$ in the vicinity of each singularity
(via each cos term in the sum),
while allowing it
to move up to $\dmax$ from the boundary in the absence of nearby singularities.
Given the curve function $y(\chi)$
a set of $N$ real values $\{\chi_j\}\in[0,2\pi]$
are then chosen such that their local spacing is proportional to $y(\chi)$.
\footnote{%
In practice this can easily be done numerically by solving the ODE
$u'(\chi) = 1/y(\chi)$ then spline
fitting to construct an approximate inverse function for $u$.}
The MFS charge points are then given by $\by_j = Z(\chi_j - iy(\chi_j))$.

We sketch our motivation for the above algorithm.
The curve \eqref{eq:adapt} can be viewed as defining a
$z$-plane curve which under some
new parametrization (annular conformal map)
is the image of a line $\text{Im}~s = C < 0$.
However, because of the curvature behavior near singularities
induced by the cos terms in \eqref{eq:adapt},
the $\tau$ for the new parametrization is large, of order 1.
One may see intuitively that it is possible to choose a parametrization
where any given singularity $\sigma$ has its imaginary part $\tau_\sigma$
pushed to infinity by choosing the conformal map corresponding to the
$\pO$ held at electrostatic potential zero while a point charge is
placed at the singularity;
this analogy informed our choice of curve.
Finally, the choice of charge point spacing approximates
the result of choosing constant spacing on the line Im$s = $ const
under the new parametrization.
This algorithm has been developed using intuition arising from electrostatics;
we do not claim that it is the best possible choice. But we obtain
very good results with it in our numerical experiments.
The parameters we give above seem to work well in a wide variety of
randomly-generated analytic curves
(the one failure mode which the algorithm does not yet
guard against is self-intersection of the resulting $\Gamma$).

In Fig.~\ref{f:splane} we illustrate
the performance of this method on a generalization of
the crescent domain given by
\be
\mbox{Generalized crescent:}\quad
Z_{GC}(s)=e^{is}-\frac{0.1}{e^{is}+a_5}-\frac{0.07+0.02i}{e^{is}+a_6}+
\frac{0.2}{e^{is}+a_7}
\label{eq:gencres}
\ee
with $a_5 = 0.9$, $a_6=-0.8-0.2i$, $a_7=-0.2+0.5i$.
The Schwarz function has three exterior pole-type
singularities, shown by '*' symbols,
with $s$-plane locations
$1/\overline{a_5}$, $1/\overline{a_6}$, and $1/\overline{a_7}$.
The two rows of subfigures contrast the effect of
$\dmax$ given by \eqref{eq:dmax} at low vs high wavenumbers.
For low wavenumber $\dmax$ is large
and the contribution of the first term in \eqref{eq:adapt}
is small. This enables the curve to reach large negative Im $s$, hence
large distances from $\pO$ and a large point spacing,
in regions away from singularities.
We observe rapid convergence, reaching $t=10^{-14}$ by $N=280$,
and no exponential coefficient growth.
As $k$ increases, $\dmax$ drops and the first term in \eqref{eq:adapt}
starts to become significant (see dotted line in Fig.~\ref{f:splane}d),
and has the effect of bringing the curve closer
to the domain. As $k\to\infty$ this term dominates and the
curve tends to a constant distance of about 4 wavelengths
from $\pO$ everywhere on the curve.
The bottom row in the figure shows the case $k=100$.
Once $N=400$ (about 3.0 degrees of freedom per wavelength
on the boundary),
convergence is rapid, reaching $10^{-10}$ at $N=525$.
Since the coefficient norm $|\bmal| \approx 10^5$ convergence halts here.
The resulting solution $u$ is shown in Fig.~\ref{f:gencresug}.
\footnote{Computation times are quoted for a single core of a 2 GHz
Intel Core Duo laptop CPU, with 2 GB RAM, running GNU/Linux, coded in MATLAB.
For evaluation of $u$ the sum \eqref{eq:mfs} was performed naively.
The MATLAB Hankel function routine is also by no means optimal,
requiring on average 2.3 microseconds per evaluation.}
Values outside $\Omega$ have been included to highlight the manner in which
the MFS points generate the field (very large coefficients are easily noticed
due to the dark, highly oscillatory bands in these parts of $\Gamma$).

Finally we test a different shape at higher wavenumber in
Fig.~\ref{f:pentafoilug}. This domain is challenging since it has,
very close to the boundary,
five exterior branch-type
singularities in the Schwarz function. We also choose non-constant boundary
data which itself has a singularity (however since it is outside $\Gamma$ there
is no need to include its contribution in the curve formula \eqref{eq:adapt}).
There are about 165 wavelengths across the domain diameter.
An error norm of order $10^{-11}$ is reached using an $N$ corresponding to
3.5 basis functions per wavelength on the perimeter.
(since due to resonance the interior $u$ values are of order $10^2$,
but $d$ as defined below \eqref{eq:mp} is of order $10^{-5}$;
these combine to give around 8 digits of relative accuracy in $u$).
We note by contrast
that for boundary element and boundary integral formulations
it is commonly stated that 10 degrees of freedom per wavelength are required
for high accuracy.
Away from singularities (or in domains with more distant singularities)
we find barely more than 2 basis functions per wavelength
are sufficient in the high $k$ limit,
similar to what was found for the disc in Section~\ref{sec:est}.
We postpone further study of this limit, and the choice of
$\dmax$ at high wavenumber, to future work.

\section{Conclusions}


The Method of Fundamental Solutions
is a powerful tool for solving the Helmholtz BVP,
but, as we have demonstrated, the achievable accuracy is limited
by the size of the
coefficient norm $|\bmal|$.
Therefore we have analysed, for the first time,
the growth of $|\bmal|$ with basis size $N$
as one converges towards Helmholtz solutions
in the disc and other analytic domains.
We emphasize that it is not the growth in the condition number
that we are studying (since this always grows exponentially with $N$),
rather the growth in the norm of the least-squares coefficient vector.
In the disc we have theorems on convergence rate (Thm.~\ref{thm:t})
and coefficient growth (Thms.~\ref{thm:n} and \ref{thm:a}),
and for analytic domains we have corresponding conjectures
(Conjs.~\ref{cnj:t} and \ref{cnj:agen})
supported by numerical experiments in many domain shapes.
These show that the success (numerical stability and hence high accuracy)
of the MFS
relies on a choice of charge curve which does not enclose any
singularities of the analytic continuation of the solution $u$.
These singularities are associated either with the analytic
continuation of the boundary data, or with the Schwarz function of the
domain.

The conclusions for optimal choice of MFS charge points are as follows.
For the unit disc, with concentric equally-spaced charge points,
a radius between $\sqrt{\rho}$ and $\rho$ is optimal (Remark~\ref{rmk:bestR}),
$\rho$ being the radius of the nearest singularity in boundary data.
For general analytic domains, charge points placed on a curve
which adapts to the singularity locations have been shown to perform
very well (and vastly outperform charge points located using equipotential
lines of the exterior conformal map).
Thus knowledge of singularities in the Schwarz function and the boundary
data are essential to achieve the best performance on analytic domains.

Our experiments show that MFS is highly competitive with boundary
integral equations, both in terms of basis size $N$ and overall simplicity.
At high wavenumber we show that in the disc
asymptotically 2 basis functions per wavelength on the boundary
are needed, and that in more complicated nonconvex domains
with nearby singularities this need only increase to about 3.5
to achieve boundary error norms close to machine precision.
In practice this unusually small $N$ results in rapid solution of the basis
coefficients, even at high wavenumber.
However, as with boundary integral methods,
the CPU time to evaluate the interior solution at roughly 10 grid points
per wavelength is much larger, especially
using a naive implementation of the sum \eqref{eq:mfs}
and MATLAB's Hankel function routine.
Replacing this evaluation of $u$ with a Fast Multipole (FMM) summation
would be a natural next step and is expected to result in a large
speedup at high wavenumbers.

We expect our findings on coefficient growth rates, and the new
adaptive charge curve algorithm,
to be easily extendable to the exterior Helmholtz scattering
problem, for which MFS has shown promise in the engineering community
\cite{Ky96,Kark01}.






\section*{Acknowledgments}
AHB is supported by the National Science Foundation under grant
DMS-0507614. TB is supported by Engineering and Physical Sciences
Research Council grant EP/D079403/1. TB wishes to thank the DMV, the
DFG, and the Shapiro Visitor Program,
for a joint travel grant in January 2007 to visit Dartmouth.
This work
benefitted from important discussions with Lehel Banjai, Leslie Greengard,
David Karkashadze, Fridon Shubitizde, and Nick Trefethen.

\appendix

\section{Bessel function asymptotics}
\label{app:bes}

We take the
standard Taylor series \cite{a+s}
\be
J_m(z) = \left(\frac{z}{2}\right)^m \sum_{k=0}^\infty
\frac{(-z^2/4)^k}{k!(m+k)!}
\ee
and in the large-$m$ limit we may approximate $(m+k)! \approx m! m^k$,
then recognize the power series for the exponential, giving
\be
J_m(z) \sim \frac{1}{m!}\left(\frac{z}{2}\right)^m e^{-z^2/4m},
\qquad m\to\infty~.
\label{eq:jasym}
\ee
Similarly the standard series
\be
Y_m(z) = -\frac{1}{\pi}\left(\frac{z}{2}\right)^{-m}
\sum_{k=0}^{m-1} \frac{(m-k-1)!(z^2/4)^k}{k!} +
\frac{2}{\pi}\ln(z/2) J_m(z) + O(z^m)
\ee
with $(m-k-1)!\approx m!/m^k$ gives
\be
Y_m(z) \sim -\frac{m!}{\pi}\left(\frac{z}{2}\right)^{-m} e^{z^2/4m},
\qquad m\to\infty~.
\ee
Neither of these asymptotic forms are given in \cite{a+s},
however \eqref{eq:jasym} has been recently noted in the physics community
\cite{landry}.
Combining these two in \eqref{eq:sm}, using the reflection
formulae, and ignoring the lower-order $J$ contribution to the Hankel function,
gives \eqref{eq:smser}.

Bessel's equation $u''+u'/r+(1-m^2/r^2)u=0$ with the Liouville
transformation $w=r^{1/2}u$ then changing variable to
$x=r/a$, with $a^2=m^2-\frac{1}{4}$, gives the ODE
\be
\frac{d^2w}{dx^2} + a^2\left(1-\frac{1}{x^2}\right)w = 0.
\ee
The WKBJ (or Liouville-Green) asymptotic approximation (Ch. 9.3 of \cite{m+f2})
for large parameter $a$ is then
\be
w(x) \sim
\spl{
(x^{-2}-1)^{-1/4} \left( A e^{a\int_x^1 \sqrt{x^{-2}-1}\, dx}
+ Be^{-a\int_x^1 \sqrt{x^{-2}-1}\, dx}\right),
& x<1 \mbox{ (evanescent)}
}
{
(1-x^{-2})^{-1/4} \left( Ce^{i a\int_1^x \sqrt{1-x^{-2}}\, dx}
+ De^{-i a\int_1^x \sqrt{1-x^{-2}}\, dx}\right),
& x>1 \mbox{ (oscillatory)}
}
\label{eq:wbkj}
\ee
where $A,B,C,D\in\mathbb{C}$ are constants.
Note that the integral in the evanescent region can be performed analytically
and is $-I_1(x)$ as defined by \eqref{eq:Iax}; the integral in the oscillatory
region is not needed since amplitude not phase is of interest.
Since the solution $w$ is continuous through
the turning point $x=1$ (even though
\eqref{eq:wbkj} breaks down), there exist connection formulae relating
the constants:
\be
C = e^{i\pi/4} A + e^{-i\pi/4} B, \qquad D = e^{i\pi/4} B.
\ee
They can be found by comparing WKBJ to large-argument asymptotics
of the Airy functions Ai and Bi on either side of the turning point
(\eg comparing 10.4.59 with 10.4.60, and 10.4.63 with 10.4.64, in \cite{a+s},
or using 9.3.91,92 of \cite{m+f2},
or the more rigorous presentation of the Gans-Jeffreys
formulae in Ch. 11 of \cite{olver}).
If $A=0$
($w$ decaying as $x$ decreases in the evanescent region)
then $|C|=|D|=|B|$ giving an amplitude of $2|B|/(1-x^{-2})^{1/4}$
in the oscillatory region.
When transformed back such a solution $u(r)$ corresponds to the
$J_m(r)$ Bessel function.
We match the asymptotic amplitude $\sqrt{2/\pi r}$ at large argument
(see 9.2.1 of \cite{a+s}) to fix $B=\frac{1}{2}$ for all $m$.
Hence the Bessel function has typical size
\be
|J_m(r)| \approx \spl{\frac{1}{2}(a^2-r^2)^{-1/4}e^{I_a(r)},&r<a}
{(a^2-r^2)^{-1/4},&r>a.}
\label{eq:Jsize}
\ee
Note that an amplitude is implied here in the oscillatory region $r>a$.
Note also that $a$ is defined above, and $I_a(r)<0$ for $r<a$.
Similarly matching the $H_m^{(1)}(r)$ Hankel function
large argument asymptotic gives $|C|=1$, $D=0$, so $|A|=1$ (which dominates),
thus typical size
\be
|H_m^{(1)}(r)| \approx \spl{(a^2-r^2)^{-1/4}e^{-I_a(r)},&r<a}
{(a^2-r^2)^{-1/4},&r>a.}
\label{eq:Hsize}
\ee
These formulae have been checked against numerical evaluations
of Bessel functions and accurately predict amplitudes or evanescent
magnitudes everywhere apart from very close to the turning point
$r=a$ where they have a weak algebraic singularity, but still
provide an upper bound.
Substituting \eqref{eq:Jsize} and \eqref{eq:Hsize} into \eqref{eq:sm}
gives the desired \eqref{eq:unif}.



\bibliography{bibfile}

\begin{thebibliography}{10}
\expandafter\ifx\csname url\endcsname\relax
  \def\url#1{\texttt{#1}}\fi
\expandafter\ifx\csname urlprefix\endcsname\relax\def\urlprefix{URL }\fi

\bibitem{a+s}
M.~Abramowitz, I.~A. Stegun, Handbook of Mathematical Functions with Formulas,
  Graphs, and Mathematical Tables, tenth edition ed., Dover, New York, 1964.

\bibitem{que}
A.~H. Barnett, Asymptotic rate of quantum ergodicity in chaotic {E}uclidean
  billiards, Comm. Pure Appl. Math. 59 (2006) 1457--88.

\bibitem{Bo85}
A.~Bogomolny, Fundamental solutions method for elliptic boundary value
  problems, SIAM Journal on Numerical Analysis 22~(4) (1985) 644--669.

\bibitem{coltonkress}
D.~Colton, R.~Kress, Inverse acoustic and electromagnetic scattering theory,
  vol.~93 of Applied Mathematical Sciences, 2nd ed., Springer-Verlag, Berlin,
  1998.

\bibitem{CoHi53}
R.~Courant, D.~Hilbert, Methods of mathematical physics. {V}ol. {I},
  Interscience Publishers, Inc., New York, N.Y., 1953.

\bibitem{Da74}
P.~J. Davis, The {S}chwarz function and its applications, The Mathematical
  Association of America, Buffalo, N. Y., 1974, the Carus Mathematical
  Monographs, No. 17.

\bibitem{En96}
R.~Ennenbach, H.~Niemeyer, The inclusion of {D}irichlet eigenvalues with
  singularity functions, Z. Angew. Math. Mech. 76~(7) (1996) 377--383.

\bibitem{FaKa98}
G.~Fairweather, A.~Karageorghis, The method of fundamental solutions for
  elliptic boundary value problems, Adv. Comput. Math. 9~(1-2) (1998) 69--95,
  numerical treatment of boundary integral equations.

\bibitem{Ga53}
P.~R. Garabedian, Applications of analytic continuation to the solution of
  boundary value problems, J. Rational Mech. Anal. 3 (1954) 383--393.

\bibitem{Ka01}
A.~Karageorghis, The method of fundamental solutions for the calculation of the
  eigenvalues of the {H}elmholtz equation, Appl. Math. Lett. 14~(7) (2001)
  837--842.

\bibitem{Kark01}
D.~Karkashadze, On status of main singularities in {3D} scattering problems,
  in: Proceedings of VIth International Seminar/Workshop on Direct and Inverse
  Problems of Electromagnetic and Acoustic Wave Theory (DIPED), Lviv, Ukraine,
  2001.

\bibitem{Ka89}
M.~Katsurada, A mathematical study of the charge simulation method. {II}, J.
  Fac. Sci. Univ. Tokyo Sect. IA Math. 36~(1) (1989) 135--162.

\bibitem{Ka90}
M.~Katsurada, Asymptotic error analysis of the charge simulation method in a
  {J}ordan region with an analytic boundary, J. Fac. Sci. Univ. Tokyo Sect. IA
  Math. 37~(3) (1990) 635--657.

\bibitem{Ka94}
M.~Katsurada, Charge simulation method using exterior mapping functions, Japan
  J. Indust. Appl. Math. 11~(1) (1994) 47--61.

\bibitem{Ka88}
M.~Katsurada, H.~Okamoto, A mathematical study of the charge simulation method.
  {I}, J. Fac. Sci. Univ. Tokyo Sect. IA Math. 35~(3) (1988) 507--518.

\bibitem{Ka96}
M.~Katsurada, H.~Okamoto, The collocation points of the fundamental solution
  method for the potential problem, Comput. Math. Appl. 31~(1) (1996) 123--137.

\bibitem{Ki88}
T.~Kitagawa, On the numerical stability of the method of fundamental solution
  applied to the {D}irichlet problem, Japan J. Appl. Math. 5~(1) (1988)
  123--133.

\bibitem{Ki91}
T.~Kitagawa, Asymptotic stability of the fundamental solution method, in:
  Proceedings of the International Symposium on Computational Mathematics
  (Matsuyama, 1990), vol.~38, 1991.

\bibitem{Ku78}
J.~R. Kuttler, V.~G. Sigillito, Bounding eigenvalues of elliptic operators,
  SIAM Journal on Mathematical Analysis 9~(4) (1978) 768--773.

\bibitem{Ky85}
A.~G. Kyurkchan, The method of auxiliary currents and sources in wave
  diffraction problems, Soviet J. Comm. Tech. Electron. 30 (1985) 48--58,
  translated from Radiotekhn. i \`Elektron. 29 (1984), no. 11, 2129--2139
  (Russian).

\bibitem{Ky96}
A.~G. Kyurkchan, B.~Y. Sternin, V.~E. Shatalov, Singularities of continuation
  of wave fields, Physics - Uspekhi 12 (1996) 1221--1242.

\bibitem{landry}
B.~Landry, E.~J. Heller, Statistical properties of many particle
  eigenfunctions, J. Phys. A 40 (2007) 9259--74.

\bibitem{Le59}
H.~Lewy, On the reflection laws of second order differential equations in two
  independent variables, Bull. Amer. Math. Soc. 65 (1959) 37--58.

\bibitem{Mi80}
R.~F. Millar, The analytic continuation of solutions to elliptic boundary value
  problems in two independent variables, J. Math. Anal. Appl. 76~(2) (1980)
  498--515.

\bibitem{Mi86}
R.~F. Millar, Singularities and the {R}ayleigh hypothesis for solutions to the
  {H}elmholtz equation, IMA J. Appl. Math. 37~(2) (1986) 155--171.

\bibitem{m+f2}
P.~Morse, H.~Feshbach, Methods of theoretical physics, volume 2, McGraw-Hill,
  1953.

\bibitem{olver}
F.~W.~J. Olver, Asymptotics and special functions, Academic Press, New York,
  1974.

\bibitem{Smyrl}
Y.-S. Smyrlis, A.~Karageorghis, Numerical analysis of the {MFS} for certain
  harmonic problems, M2AN Math. Model. Numer. Anal. 38~(3) (2004) 495--517.

\end{thebibliography}

%
%
%
%
%

\end{document}